\documentclass[aop,citesort,MSNbibl,dvips]{arximspdf}
\usepackage{graphics}
\doi{10.1214/09-AOP524}
\volume{38}
\issue{5}
\pubyear{2010}
\firstpage{1783}
\lastpage{1816}

\makeatletter

\newcommand{\llbracket}{[\![}
\newcommand{\rrbracket}{]\!]}

\newproclaim{defin}{Definition}
\newtheorem{theo}[defin]{Theorem}
\newproclaim{nota}[defin]{Notation}
\newtheorem{cor}[defin]{Corollary}
\newtheorem{prop}[defin]{Proposition}
\newtheorem{lem}[defin]{Lemma}

\makeatother

\begin{document}
\begin{frontmatter}

\title{Asymptotics of one-dimensional forest fire processes}
\runtitle{Forest fire processes}

\begin{aug}
\author[A]{\fnms{Xavier} \snm{Bressaud}\corref{}\ead[label=e1]{bressaud@math.univ-toulouse.fr}} and
\author[B]{\fnms{Nicolas} \snm{Fournier}\thanksref{t1}\ead[label=e2]{nicolas.fournier@univ-paris12.fr}}
\runauthor{X. Bressaud and N. Fournier}
\affiliation{Universit\'e Paul Sabatier and Universit\'e Paris-Est}
\address[A]{Institut de Math\'ematiques de Toulouse\\
Universit\'e Paul Sabatier\\
F-31062 Toulouse Cedex 9\\
France\\
\printead{e1}} 
\address[B]{LAMA, Facult\'e de Sciences et Technologie\\
Universit\'e Paris-Est\\
61, Avenue du G\'en\'eral de Gaulle\\
94010 Cr\'eteil Cedex\\
France\\
\printead{e2}}
\end{aug}

\thankstext{t1}{Supported by  Agence Nationale de la Recherche
Grant ANR-08-BLAN-0220-01.}

\received{\smonth{12} \syear{2008}}
\revised{\smonth{9} \syear{2009}}

%
\begin{abstract}
We consider the so-called one-dimensional forest fire process.
At each site of $\mathbb{Z}$, a tree appears at rate $1$.
At each site of $\mathbb{Z}$, a fire starts at rate ${\lambda}>0$,
immediately destroying the whole corresponding connected component of trees.
We show that when ${\lambda}$ is made to tend to $0$ with an appropriate
normalization, the forest fire process tends to a uniquely
defined process, the dynamics of which we precisely describe.
The normalization consists of accelerating time by
a factor $\log(1/{\lambda})$ and of compressing space by a factor
${\lambda}\log(1/{\lambda})$.
The limit process is quite simple: it can be built using a graphical
construction and can be perfectly simulated. Finally, we derive some
asymptotic estimates (when ${\lambda}\to0$)
for the cluster-size distribution of the forest fire process.
\end{abstract}

%

\end{frontmatter}

\section{Introduction and main results}

\subsection{The model}

Consider two independent families of independent Poisson processes,
$N=(N_t(i))_{t\geq0,i\in\mathbb{Z}}$ and
$M^{\lambda}=(M^{\lambda}_t(i))_{t\geq0,i\in\mathbb{Z}}$, with
respective rates~$1$ and ${\lambda}>0$. Define
${\mathcal F}_t^{N,M^{\lambda}}:=\sigma(N_s(i),M^{\lambda}_s(i),
s\leq t, i \in\mathbb{Z})$.
For $a,b \in\mathbb{Z}$ with \mbox{$a\leq b$}, we set $\llbracket a, b
\rrbracket
= \{a,\ldots,b\}$.
\begin{defin}\label{dflaffp}
Consider a ${\{0,1\}^\mathbb{Z}}$-valued $({\mathcal
F}_t^{N,M^{\lambda}})_{t\geq0}$-adapted
process\break $(\eta^{\lambda}_t)_{t\geq0}$ such that $(\eta^{\lambda
}_t(i))_{t\geq0}$
is a.s. c\`adl\`ag for all $i\in\mathbb{Z}$.

We say that $(\eta^{\lambda}_t)_{t\geq0}$
is a ${\lambda}$-FFP (forest fire process)
if a.s., for all $t\geq0$ and all $i\in\mathbb{Z}$,
\[
\eta^{\lambda}_t(i)={\int_0^t}{\mathbf{1}}_{\{\eta^{\lambda
}_{{s-}}(i)=0\}} \,dN_s(i)
- \sum_{k\in\mathbb{Z}}{\int_0^t}{\mathbf{1}}_{\{k\in C^{\lambda
}_{{s-}}(i)\}} \,dM^{\lambda}_s(k),
\]
where $C^{\lambda}_s(i)=\varnothing$ if $\eta_t^{\lambda}(i)=0$, while
$C^{\lambda}_s(i)=\llbracket l_s^{\lambda}(i),r_s^{\lambda
}(i)\rrbracket$ if $\eta
^{\lambda}
_s(i)=1$, with
\[
l_s^{\lambda}(i)=\sup\{k< i; \eta_s^{\lambda}(k)=0\}+1
\quad\mbox{and}\quad
r_s^{\lambda}(i)=\inf\{k > i; \eta_s^{\lambda}(k)=0\}-1.
\]
\end{defin}

Formally, we say that $\eta_t^{\lambda}(i)=0$ if there is no tree at
site $i$
at time $t$ and $\eta_t^{\lambda}(i)=1$ otherwise. $C_t^{\lambda
}(i)$ stands for
the connected component of occupied sites around $i$ at time $t$.
Thus, the forest fire process starts
from an empty initial configuration, trees appear on vacant sites at
rate $1$ (according to $N$) and a fire starts on each site at rate
${\lambda}>0$ (according to $M^{\lambda}$), immediately burning the
corresponding
connected component of occupied sites.

This process can be shown to exist and to be unique (for almost every
realization of $N,M^{\lambda}$) by using a \textit{graphical construction}.
Indeed, to build the process until a given time $T>0$,
it suffices to work between sites $i$ which are vacant until time $T$
[because $N_T(i)=0$]. Interaction cannot cross such sites. Since
such sites are a.s. infinitely many, this allows us to handle
a graphical construction.
We refer to Van den Berg and Jarai \cite{vdbj} (see also
Liggett \cite{l}) for many examples of graphical constructions.
It should be pointed out that this construction only works in
dimension~$1$.

\subsection{Motivation and references}\label{bibi}

The study of self-organized critical (SOC) systems has
become rather popular in physics since the end of the 1980s.
SOC systems are simple models which are supposed to shed light on
temporal and spatial randomness observed in a variety of natural
phenomena showing \textit{long-range correlations}, like sand piles, avalanches,
earthquakes, stock market crashes, forest fires, shapes of mountains,
clouds, etc. Roughly,
the idea, which appears in Bak, Tang and Wiesenfeld \cite{btw1}
with regard to sand piles, is that of
systems \textit{growing} toward a \textit{critical state}
and relaxing through \textit{catastrophic}
events (avalanches, crashes, fires, etc.).
The most classical model is the sand pile model introduced
in 1987 in~\cite{btw1}, but many variants
or related models have been proposed and studied more or less rigorously,
describing earthquakes (Olami, Feder and Christensen \cite{ofc}) or
forest fires (Henley \cite{h}, Drossel and Schwabl \cite{ds}).
For surveys on the subject, see
Bak, Tang and Wiesenfeld \cite{btw1,btw2}, Jensen
\cite{j} and the references therein.

From the point of view of SOC systems, the forest fire model is
interesting in the
asymptotic regime $\lambda
\to0$. Indeed, fires are less frequent, but when they occur, destroyed
clusters may be huge. This model has been the subject of many numerical
and heuristic studies; see Drossel, Clar and Schwabl \cite{dcs}
and Grassberger \cite{g} for references. However,
there are few rigorous results. Even existence
of the (time-dependent) process for a multidimensional
lattice and given $\lambda>0$
has been proven only recently \cite{du1,du2} and uniqueness is
known to hold only for ${\lambda}$ large enough.
The existence and uniqueness of an invariant distribution
(as well as other qualitative properties),
even in dimension $1$,
have been proven only recently in \cite{bf} for ${\lambda}=1$. These last
results can probably be extended to the case where $\lambda\geq1$,
but the
method in \cite{bf}
completely breaks down for small values of ${\lambda}$.

The asymptotic behavior of the ${\lambda}$-FFP as ${\lambda}\to0$
has been studied
numerically and heuristically \cite{ds,dcs,d,g}.
To our knowledge, the only mathematically rigorous results are the following.

(a) Van den\vspace*{2pt} Berg and Jarai \cite{vdbj} have proven that
for $t\geq3$, ${\mathbb{P}}[\eta_{t\log(1/{\lambda})}^{\lambda
}(0)=0]\simeq1/\log
(1/{\lambda})$,
thus giving some idea of the density of vacant sites. This result
was conjectured by Drossel, Clar and Schwabl \cite{dcs}.

(b) Van den Berg and Brouwer \cite{vdbb} have obtained some
results in the two-dimensional case concerning the behavior
of clusters near the \textit{critical time}. However, these results
are not completely rigorous since they are based on a percolation-like
assumption, which is not rigorously proved.

(c) Brouwer and Pennanen \cite{bp} have proven the existence of an
invariant distribution for each fixed ${\lambda}>0$, as well as a precise
version of the following
estimate which extends (a): for ${\lambda}\in(0,1)$, at equilibrium,
${\mathbb{P}}[\#(C^{\lambda}(0))=x ] \simeq c /[x \log(1/{\lambda})]$
for $x \in\{ 1,\ldots,(1/{\lambda})^{1/3}\}$. It was conjectured
in \cite{dcs} that this actually
holds for $x \in\{1,\ldots,1/({\lambda}\log(1/{\lambda}))\}$,
but this was rejected in \cite{vdbj}.

In this paper, we rigorously derive a limit theorem which shows that
the ${\lambda}$-FFP converges, under rescaling, to some limit forest
fire process
(LFFP). We precisely describe the dynamics of the LFFP and
show that it is quite simple: in particular, it is unique,
can be built by using a \textit{graphical construction} and can thus
be \textit{perfectly} simulated. Our result
allows us to prove a very weak version of (c)
for $x\in\{1,\ldots,(1/{\lambda})^{1-{\varepsilon}}\}$, for any
${\varepsilon}>0$;
see Corollary \ref{coco} below.

\subsection{Notation}\label{ssnota}

We denote by $\#(I)$ the number of elements of a set $I$.

For $a,b \in\mathbb{Z}$, with $a\leq b$, we set $\llbracket a, b
\rrbracket
= \{a,\ldots,b\}\subset\mathbb{Z}$.

For $I=\llbracket a,b\rrbracket\subset\mathbb{Z}$
and
$\alpha>0$, we will set $\alpha I := [\alpha a, \alpha b]\subset
{\mathbb{R}}$.
For $\alpha>0$, we naturally adopt the convention that
$\alpha\varnothing=\varnothing$.

For $J=[a,b]$, an interval of ${\mathbb{R}}$,
$|J|=b-a$ stands for the length of $J$ and for $\alpha>0$, we set
$\alpha J = [\alpha a, \alpha b]$.

For $x\in{\mathbb{R}}$, $\lfloor x \rfloor$ stands for the integer
part of $x$.

\subsection{Heuristic scales and relevant quantities}\label{hscales}

Our aim is to find some time scale for which tree clusters experience
approximately
one fire per unit of time. However, for ${\lambda}$ very small,
clusters will
be very large immediately before they burn. We must thus also rescale space,
in order that, immediately before  burning, clusters have a size
of order $1$.

\subsubsection*{Time scale} Consider the cluster $C_t^{\lambda}(x)$ around some
site $x$
at time $t$.
It is quite clear that for ${\lambda}>0$ very small and $t$ not too large,
one can neglect fires
so that, roughly, each site is occupied with probability $1-e^{-t}$
and, thus, $C^{\lambda}_t(x)\simeq\llbracket x-X,x+Y\rrbracket$,
where $X,Y$
are geometric
random variables with parameter $1-e^{-t}$.
As a consequence, $\#(C_t^{\lambda}(x))\simeq e^{t}$ for $t$ not too large.
On the other hand, the cluster $C_t^{\lambda}(x)$ burns at rate
${\lambda}\#
(C_t^{\lambda}(x))$
(at time $t$) so that we decide to accelerate time by
a factor $\log(1/{\lambda})$. In this way, ${\lambda}\#(C_{\log
(1/{\lambda})}^{\lambda}
(x))\simeq1$.

\subsubsection*{Space scale}
We now rescale space in such a way that during a time interval of order
$\log(1/{\lambda})$, something like one fire starts per unit of
(space) length.
Since fires occur at rate ${\lambda}$, our space scale has to be of order
${\lambda}\log(1/{\lambda})$: this means that we will identify
$\llbracket0,\lfloor1/({\lambda}\log(1/{\lambda}))\rfloor
\rrbracket\subset
\mathbb{Z}$
with $[0,1]\subset{\mathbb{R}}$.

\subsubsection*{Rescaled clusters} We thus set, for ${\lambda}\in(0,1)$,
$t\geq0$ and $x\in{\mathbb{R}}$,
recalling Section \ref{ssnota},
%
%
\begin{equation}\label{dlambda}
D^{\lambda}_t(x):= {\lambda}\log(1/{\lambda}) C^{\lambda}_{t\log
(1/{\lambda})}\bigl(
\bigl\lfloor x /\bigl( {\lambda}\log(1/{\lambda})\bigr) \bigr\rfloor\bigr)
\subset{\mathbb{R}}.
\end{equation}

However, this creates an immediate difficulty:
recalling that $\#(C_t^{\lambda}(x)) \simeq e^t$ for $t$ not too large,
we see that for each site $x$,
$|D^{\lambda}_t(x)| \simeq{\lambda}\log(1/{\lambda}) e^{t \log
(1/{\lambda})}=
{\lambda}^{1-t} \log(1/{\lambda})$, of which the limit as ${\lambda
}\to0$ is
$0$ for $t<1$ and $+\infty$ for $t\geq1$.

For $t\geq1$, there might be fires in effect
and one hopes that this will make the possible limit of $|D^{\lambda
}_t(x)|$ finite.
However, fires can only reduce the size of clusters so that for $t<1$,
the limit of $|D^{\lambda}_t(x)|$ will really be $0$.
Thus, for a possible limit $|D(x)|$
of $|D^{\lambda}(x)|$, we should observe some paths of the following form:
$|D_t(x)|=0$ for $t<1$, $|D_t(x)|>0$ for some times $t\in(1,\tau)$,
after which it might be killed by a fire
and thus come back to $0$, at which time it remains at $0$ for
a time interval of length $1$, and so on.

This cannot be a Markov process
because $|D(x)|$ always remains at $0$
during a time interval of length exactly $1$.
We thus need to keep track of more information in order to control when
it exits from $0$.

\subsubsection*{Degree of smallness}
As was stated previously, we hope that for $t<1$,
$|D^{\lambda}_t(x)| \simeq{\lambda}^{1-t} \log(1/{\lambda})\simeq
{\lambda}^{1-t}$.
Thus, we will try to keep in mind the degree of smallness.
We will define, for
${\lambda}\in(0,1)$, $x\in{\mathbb{R}}$ and $t>0$,
%
\begin{equation}\label{zlambda}
Z^{\lambda}_t(x):= \frac{\log[1+\# ( C^{\lambda}_{t\log(1/{\lambda})}(
\lfloor x /( {\lambda}\log(1/{\lambda})) \rfloor) ) ]}
{\log(1/{\lambda})} \in[0,\infty).
\end{equation}

\subsubsection*{Final description}
We will study the ${\lambda}$-FFP via
$(D^{\lambda}_t(x),Z^{\lambda}_t(x))_{x\in{\mathbb{R}},t\geq0}$.
The main idea is that for ${\lambda}>0$ very small:

\begin{longlist}
\item if $Z^{\lambda}_t(x)=z\in(0,1)$, then $|D^{\lambda}_t(x)|\simeq
0$ and the (rescaled)
cluster containing
$x$ is microscopic, but we control its smallness, in the sense that
$|D^{\lambda}_t(x)|\simeq{\lambda}^{1-z}$ (in a very unprecise way);

\item if $Z^{\lambda}_t(x)=1$ [we will show below that $Z^{\lambda
}_t(x)$ will never
exceed $1$ in the limit ${\lambda}\to0$], then
the (rescaled) cluster containing $x$ is automatically macroscopic and
has a length
equal to $|D^{\lambda}_t(x)|\in(0,\infty)$.
\end{longlist}

\subsection{The limit process}\label{sslffp}

We now describe the limit process. We want this process to be Markov
and this forces us to add some variables.

We consider a Poisson measure $M(dt,dx)$ on $[0,\infty) \times
{\mathbb{R}}$ with
intensity measure $dt \,dx$. Again, we define
${\mathcal F}_t^M=\sigma(M(A),A\in{\mathcal B}([0,t] \times{\mathbb{R}}))$.
We also define ${\mathcal I}:=\{[a,b], a \leq b\}$,
the set of all closed finite intervals of ${\mathbb{R}}$.
\begin{defin}\label{dflffp}
A $({\mathcal F}_t^M)_{t\geq0}$-adapted process
$(Z_t(x),D_t(x),H_t(x))_{t\geq0,x\in{\mathbb{R}}}$ with
values in ${\mathbb{R}}_+\times{\mathcal I}\times{\mathbb{R}}_+$ is
a limit forest fire
process (LFFP)
if a.s., for all $t\geq0$ and all $x \in{\mathbb{R}}$,
%
%
\begin{equation}\label{eqlffp}
\cases{
\displaystyle Z_t(x)= {\int_0^t}{\mathbf{1}}_{\{Z_s (x) < 1\}}
\,ds -
{\int_0^t}\int_{\mathbb{R}}{\mathbf{1}}_{\{ Z_{{s-}}(x)=1,y \in
D_{{{s-}}}(x)\}
}M(ds,dy),\cr
\displaystyle H_t(x)=
{\int_0^t}Z_{{s-}}(x){\mathbf{1}}_{\{Z_{{s-}}(x)<1\}} M(ds\times\{x\})
- {\int_0^t}{\mathbf{1}}_{\{H_s (x) > 0 \}}\,ds,}
\end{equation}
where $D_t(x) = [L_t(x),R_t(x)]$ with
\begin{eqnarray*}
L_t(x) &=& \sup\{ y\leq x; Z_t(y)<1 \mbox{ or } H_t(y)>0 \},\\
R_t(x) &=& \inf\{ y\geq x; Z_t(y)<1 \mbox{ or } H_t(y)>0 \}.
\end{eqnarray*}
\end{defin}

A typical path of the finite box version of the LFFP
(see Section \ref{ex}) is drawn and commented on in Figure \ref{figure2}
and a simulation algorithm is explained in the proof of Proposition
\ref{eufini}.

Let us explain the dynamics of this process. We consider $T>0$ fixed
and set ${\mathcal B}_T= \{x \in{\mathbb{R}}; M([0,T]\times\{x\})>
0\}$.
For each $t\geq0$ and $x\in{\mathbb{R}}$, $D_t(x)$ stands for the
occupied cluster containing $x$. We call this cluster \textit{microscopic}
if $D_t(x)=\{x\}$. We also have $D_t(x)=D_t(y)$ for all
$y$ in the interior of $D_t(x)$: if $D_t(x)=[a,b]$, then $D_t(y)=[a,b]$
for all $y\in(a,b)$.

1. \textit{Initial condition.}
We have $Z_0(x)=H_0(x)=0$
and $D_0(x)=\{x\}$ for all $x\in{\mathbb{R}}$.

2. \textit{Occupation of vacant zones.}
Here, we consider $x\in{\mathbb{R}}\setminus{\mathcal B}_T$. We
then have
$H_t(x)=0$ for all $t\in[0,T]$.
If $Z_t(x)<1$, then $D_t(x)=\{x\}$
and $Z_t(x)$ stands for the \textit{degree of
smallness} of the cluster containing $x$. Then
$Z_t(x)$ grows linearly until it reaches $1$, as described by the first term
on the right-hand side of the first equation in (\ref{eqlffp}). If
$Z_t(x)=1$, then
the cluster containing $x$ is macroscopic and is described
by $D_t(x)$.

3. \textit{Microscopic fires.} Here, we assume that $x\in{\mathcal B}_T$
and that the
corresponding mark
of $M$ happens at some time $t$ where $z:=Z_{{t-}}(x)<1$. In such a case,
the cluster containing $x$ is microscopic. We then set $H_t(x)=Z_{{t-}}(x)$,
as described by the first term on the right-hand side of the second
equation of
(\ref{eqlffp}), and we leave the value of $Z_t(x)$ unchanged.
We then let $H_s(x)$ decrease linearly until it reaches~$0$;
see the second term on the right-hand side of the second equation in
(\ref{eqlffp}).
At all times where $H_s(x)>0$, that is, during $[t,t+z)$,
the site $x$ acts like a barrier (see point 5
below).

4. \textit{Macroscopic fires.}
Here, we assume that $x\in{\mathcal B}_T$ and that the corresponding mark
of $M$ happens at some time $t$ where $Z_{{t-}}(x)=1$. This means that
the cluster containing $x$ is macroscopic and thus this mark destroys
the whole component $D_{{t-}}(x)$. That is, for all $y\in D_{{t-}}(x)$,
we set
$D_t(y)=\{y\}$, $Z_t(y)=0$. This is described by the second term on the
right-hand side
of the first equation in (\ref{eqlffp}).

5. \textit{Clusters.} Finally, the definition of the clusters
$(D_t(x))_{x\in{\mathbb{R}}}$ becomes more clear: these
clusters are delimited by zones with microscopic sites [i.e., $Z_t(y)<1$]
or by sites where there has (recently) been a microscopic fire
[i.e., $H_t(y)>0$].

\subsection{Main results}

First, we must note that it is not entirely clear that the limit
process exists.
\begin{theo}\label{existe}
For any Poisson measure $M$, there a.s. exists a unique LFFP;
recall Definition \ref{dflffp}.
Furthermore, it can be constructed graphically and thus its restriction
to any finite box $[0,T]\times[-n,n]$ can be perfectly simulated.
\end{theo}

To describe the convergence of the ${\lambda}$-FFP to the LFFP, we
will need some more notation. Let ${\mathbb{D}}([0,T],E)$ denote the
space of
right-continuous and left-limited functions from the interval $[0,T]$
to a topological space $E$.
\begin{nota}\label{nocv}
(i) For two intervals $[a,b]$ and $[c,d]$, we set
$\delta([a,b],[c,d])=|a-c|+|b-d|$. We also set, by convention,
$\delta([a,b],\varnothing)=|b-a|$.

(ii) For $(x,I), (y,J)$ in ${\mathbb{D}}([0,T],
{\mathbb{R}}\times{\mathcal I}\cup\{\varnothing\})$, let
\[
\delta_T((x,I),(y,J))={\sup_{[0,T]}}|x(t)-y(t)|+ \int_0^T\delta
(I(t),J(t)) \,dt.
\]
\end{nota}

We are finally in a position to state our main result.
\begin{theo}\label{converge}
Consider, for all ${\lambda}>0$, the processes
$(Z^{\lambda}_t(x),D^{\lambda}_t(x))_{t\geq0,x\in{\mathbb{R}}}$
associated with
the ${\lambda}$-FFP;
see Definition \ref{dflaffp} and (\ref{dlambda}), (\ref{zlambda}).
Let $(Z_t(x),D_t(x)$,\break $H_t(x))_{t\geq0,x\in{\mathbb{R}}}$ be an LFFP,
as in
Definition \ref{dflffp}.

\textup{(a)} For any $T>0$ and any finite subset $\{x_1,\ldots,x_p\}\subset
{\mathbb{R}}$,
$(Z^{\lambda}_t(x_i)$,\break $D^{\lambda}_t(x_i))_{t\in[0,T],i=1,\ldots,p}$ goes
in law to
$(Z_t(x_i),D_t(x_i))_{t\in[0,T],i=1,\ldots,p}$ in ${\mathbb{D}}([0,T]$,
${\mathbb{R}}\times
{\mathcal I})^p$
as ${\lambda}$ tends to $0$. Here,
${\mathbb{D}}([0,\infty)$, ${\mathbb{R}}\times{\mathcal I})$ is
endowed with the distance
$\delta_T$;
see Notation \ref{nocv}.

\textup{(b)} For any finite subset $\{(t_1,x_1),\ldots,(t_p,x_p)\}\subset
{\mathbb{R}}_+ \times{\mathbb{R}}$,
$(Z^{\lambda}_{t_i}(x_i)$,\break $D^{\lambda}_{t_i}(x_i))_{i=1,\ldots,p}$ goes in
law to
$(Z_{t_i}(x_i),D_{t_i}(x_i))_{i=1,\ldots,p}$ in $({\mathbb{R}}\times
{\mathcal I})^p$.
\end{theo}

Observe that the process $H$ does not
appear in the limit since for each $x\in{\mathbb{R}}$, a.s., for all
$t\geq0$,
$H_t(x)=0$. [Of course, it is not the case that a.s., for all $x\in
{\mathbb{R}}$,
all $t\geq0$, $H_t(x)=0$.]
We obtain the convergence of $D^{\lambda}$ to $D$ only when
integrating in time.
We cannot hope for a Skorokhod convergence
since the limit process $D(x)$ jumps instantaneously from $\{x\}$
to some interval with positive length, while $D^{\lambda}(x)$ needs many
small jumps
(in a very short time interval) to become macroscopic.

As a matter of fact, we will obtain a convergence in probability,
using a coupling argument. Essentially, we will consider a
Poisson measure $M(dt,dx)$, as in Section \ref{sslffp}, and set, for
${\lambda}\in(0,1)$ and $i \in\mathbb{Z}$,
\[
M^{\lambda}_t(i)=M \bigl([0,t/\log(1/{\lambda})]\times
\bigl[i{\lambda}\log(1/{\lambda}),(i+1){\lambda}\log(1/{\lambda})\bigr) \bigr).
\]
Then $(M^{\lambda}_t(i))_{t\geq0, i\in\mathbb{Z}}$ is an i.i.d.
family of
Poisson processes with rate ${\lambda}$.

The i.i.d. family of Poisson processes $(N_t(i))_{t\geq0, i\in\mathbb{Z}}$
with rate
$1$ can be chosen arbitrarily, but we will decide to choose the same
family for all values of ${\lambda}\in(0,1)$.

\subsection{Heuristic arguments}

We now explain roughly the reasons why Theorem~\ref{converge}
holds.
We consider a
${\lambda}$-FFP $(\eta^{\lambda}_t)_{t\geq0}$ and
the associated process $(Z^{\lambda}_t(x)$, $D^{\lambda}_t(x))_{t\geq
0,x\in{\mathbb{R}}}$.
We assume below that ${\lambda}$ is very small.

0. \textit{Scales.} With our scales,
there are $1/({\lambda}\log(1/{\lambda}))$ sites per unit of length.
Approximately
one fire starts per unit of time per unit of length. A vacant site
becomes occupied at rate $\log(1/{\lambda})$.

1. \textit{Initial condition.}
We have, for all $x\in{\mathbb{R}}$,
$(Z^{\lambda}_0(x),D^{\lambda}_0(x))= (0,\varnothing) \simeq(0,\{x\})$.

2. \textit{Occupation of vacant zones.}
Assume that a zone $[a,b]$ (which corresponds to the zone
$\llbracket\lfloor
a/({\lambda}\log(1/{\lambda}))\rfloor,b/({\lambda}\log(1/{\lambda
}))\rfloor\rrbracket$
before rescaling)
becomes
completely vacant at some time $t$ [or $t\log(1/{\lambda})$ before rescaling]
because it has been destroyed
by a fire.

{\smallskipamount=0pt
\begin{longlist}
\item For $s\in[0,1)$, and if no fire starts on $[a,b]$ during
$[t,t+s]$, we have $D^{\lambda}_{t+s}(x) \simeq[x \pm{\lambda}^{1-s}]$
and thus $Z^{\lambda}_{t+s}(x)\simeq s$
for all $x\in[a,b]$.

Indeed, $D^{\lambda}_{t+s}(x)\simeq[x-{\lambda}\log(1/{\lambda})X,
x+{\lambda}\log(1/{\lambda}) Y]$,
where $X$ and $Y$ are geometric random variables with parameter
$1-e^{-s\log(1/{\lambda})}=1-{\lambda}^s$. This comes from the fact
that each site
of $[a,b]$ is vacant at time $t$ and becomes occupied at rate
$\log(1/{\lambda})$.

\item If no fire starts on $[a,b]$ during $[t,t+1]$, then
$Z^{\lambda}_{t+1}(x)\simeq1$ and all the sites in $[a,b]$ are occupied
(with very high probability) at time $t+1$. Indeed,
we have $(b-a)/({\lambda}\log(1/{\lambda}))$ sites and each of them
is occupied
at time $t+1$ with probability $1-e^{-\log(1/{\lambda})}=1-{\lambda}$
so that all of them are occupied with
probability $(1-{\lambda})^{(b-a)/({\lambda}\log(1/{\lambda
}))}\simeq
e^{-(b-a)/\log(1/{\lambda})}$, which goes to $1$ as ${\lambda}\to0$.
\end{longlist}}

3. \textit{Microscopic fires.}
Assume that a fire starts at some location $x$ (i.e., $\lfloor
x/({\lambda}\log(1/{\lambda}))\rfloor$ before rescaling)
at some time $t$ [or $t\log(1/{\lambda})$ before rescaling]
with $Z_{t-}^{\lambda}(x)=z\in(0,1)$.
The possible clusters on the left and
right of $x$ cannot then be connected during (approximately) $[t,t+z]$,
but they can be connected after (approximately) $t+z$. In other words,
$x$ acts like a barrier during $[t,t+z]$.

Indeed, the fire makes vacant a zone $A$ of approximate length
${\lambda}^{1-z}$ around $x$,
which thus contains approximately ${\lambda}^{1-z}/({\lambda}\log
(1/{\lambda}))
\simeq
{\lambda}^{-z}$ sites.
The probability that a fire starts again
in $A$ after $t$ is very small. Thus, using the same
computation
as in point 2(ii), we observe that ${\mathbb{P}}[A $ is
completely occupied
at time $ t+s]\simeq(1-{\lambda}^s)^{{\lambda}^{-z}}\simeq
e^{-{\lambda}^{s-z}}$. When
${\lambda}\to0$, this quantity
tends to $0$ if $s<z$ and to $1$ if $s>z$.

4. \textit{Macroscopic fires.}
Assume, now, that a fire starts at some place $x$ (i.e., $\lfloor
x/({\lambda}\log(1/{\lambda}))\rfloor$ before rescaling)
at some time $t$ [or $t\log(1/{\lambda})$ before rescaling]
and that $Z^{\lambda}_{t}(x) \simeq1$. Thus,
$D^{\lambda}_{t}(x)$ is macroscopic (i.e., its length is of order $1$
in our scales). This will thus make vacant the zone $D^{\lambda
}_t(x)$. Such a
(macroscopic) zone
needs a time of order $1$ to be completely occupied, as explained in
point 2(ii).

5. \textit{Clusters.} For $t\geq0$, $x\in{\mathbb{R}}$, the cluster
$D_t^{\lambda}(x)$
resembles $[x\pm{\lambda}^{1-z}]\simeq\{x\}$ if $Z_t^{\lambda
}(x)=z\in(0,1)$.
We then say that $x$ is microscopic.
Now, macroscopic clusters are delimited either by
microscopic zones or by
sites where there has been a microscopic fire (see point 3).

Comparing the arguments above to the rough description
of the LFFP (see Section \ref{sslffp}), our hope is that the ${\lambda}$-FFP
resembles the LFFP for ${\lambda}>0$ very small.

\subsection{Decay of correlations}

A byproduct of our result is an estimate on the decay of correlations
in the LFFP for finite times. We refer to Proposition \ref{plocla} below
for a precise statement. The main idea is that for all $T>0$,
there are constants $C_T>0$, $\alpha_T>0$ such that for all ${\lambda
}\in
(0,1)$ and
all $A>0$, the values of the
${\lambda}$-FFP inside $[-A/({\lambda}\log(1/{\lambda
})),A/({\lambda}\log(1/{\lambda}))]$
are independent of the values outside $[-2A/({\lambda}\log(1/{\lambda})),
2A/({\lambda}\log(1/{\lambda}))]$
during the time interval $[0,T\log(1/{\lambda})]$, up to a
probability smaller
that $C_Te^{-\alpha_T A}$. In other words,
for times of order $\log(1/{\lambda})$, the range of correlations
is at most of order $1/({\lambda}\log(1/{\lambda}))$.

\subsection{Cluster size distribution}

Finally, we give results on the cluster size distribution,
which are to be compared with \cite{vdbj,bp};
see Section \ref{bibi} above.
\begin{cor}\label{coco}
For each ${\lambda}>0$, consider a ${\lambda}$-FFP process $(\eta
^{\lambda}_t)_{t\geq0}$.

\begin{longlist}
\item For some $0<c<C$,
all $t\geq5/2$ and all $0\leq a < b < 1$,
\[
c (b-a) \leq\lim_{{\lambda}\to0}
{\mathbb{P}}\bigl(\#\bigl(C^{\lambda}_{t\log(1/{\lambda})}(0)\bigr) \in
[{\lambda}^{-a},{\lambda}^{-b}] \bigr) \leq C (b-a).
\]

\item
For some $0<c<C$, some $0< \kappa_1 <\kappa_2$, all
$t \geq3/2$ and all $B>0$,
\[
c e^{-\kappa_2 B} \leq\lim_{{\lambda}\to0}
{\mathbb{P}}\bigl(\#\bigl(C^{\lambda}_{t\log(1/{\lambda})}(0)\bigr)
\geq B / \bigl({\lambda}\log(1/{\lambda})\bigr) \bigr)
\leq C e^{-\kappa_1 B}.
\]
\end{longlist}
\end{cor}

Point (i) says, roughly, that for $t$ large enough
(say at equilibrium) and
for $x << 1/{\lambda}$ [say for $x \leq(1/{\lambda})^{1-{\varepsilon}}$],
choosing $a=\log(x) / \log(1/{\lambda})$ and $b=\log(x+1) / \log
(1/{\lambda})$,
we have
\begin{eqnarray*}
{\mathbb{P}}\bigl(\#(C^{\lambda}(0)) = x\bigr) &\simeq&{\mathbb{P}}\bigl(\#
(C^{\lambda}(0)) \in[x,x+1]\bigr)\simeq
{\mathbb{P}}\bigl(\#(C^{\lambda}(0)) \in[{\lambda}^{-a},{\lambda
}^{-b}]\bigr)\\
&\simeq&(b-a) \simeq\frac{1}
{x\log(1/{\lambda})}.
\end{eqnarray*}
It is thus a very weak form
of the result of \cite{bp}, but it holds for
a much wider class of $x$:
here, we allow $x \leq1/{\lambda}^{1-{\varepsilon}}$, while $x\leq
1/{\lambda}^{1/3}$
was imposed in \cite{bp}. Another advantage of our result
is that we can prove that the limit exists in (i).

%
\begin{figure}

\includegraphics{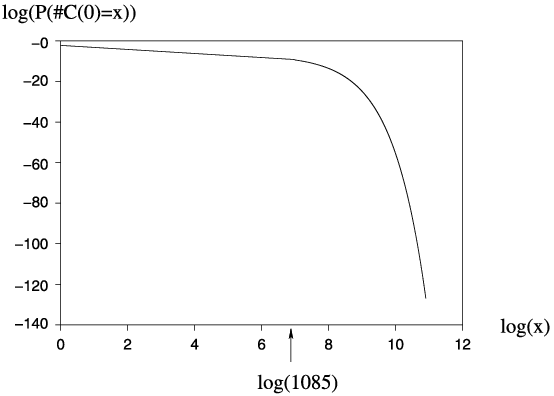}

\caption{Shape of the cluster size distribution.
Here, ${\lambda}=0.0001$ and the critical size is thus
$1/({\lambda}\log(1/{\lambda}))\simeq1085$.
We have drawn the approximate value (computed roughly just after
Corollary \protect\ref{coco}) of $\log({\mathbb{P}}(\#(C^{\lambda
}(0))=x))$ as a function
of $\log(x)$ for $x=1,\ldots,54\mbox{,}250$.
We have made the curve continuous around $x=1085$ (without justification).
The curve is linear for $x=1,\ldots,1085$ and nonlinear for
$x\geq1085$.}
\label{figure1}
\end{figure}

Point (ii) roughly describes the cluster size distribution
of macroscopic components, that is, of components of which the
size is of order $1/({\lambda}\log(1/{\lambda}))$. Here, again,
rough computations
show that for $x > {\varepsilon}/({\lambda}\log(1/{\lambda}))$ and
for $t$
large enough
(say at equilibrium),
\[
{\mathbb{P}}\bigl(\#(C^{\lambda}(0)) = x\bigr) \simeq{\lambda}\log
(1/{\lambda}) e^{-\kappa x {\lambda}\log
(1/{\lambda})}.
\]
Thus, there is clearly a phase transition near the \textit{critical size}
$1/({\lambda}\log(1/{\lambda}))$; see Figure \ref{figure1} for an
illustration.

\subsection{Organization of the paper}

The paper is organized as follows. In Section~\ref{ex}, we give the proof
of Theorem \ref{existe}. In Section \ref{locla}, we show
that, in some sense, the ${\lambda}$-FFP can be localized in a finite
box, uniformly for ${\lambda}>0$.
Section \ref{cv} is devoted to the proof of Theorem \ref{converge}.
Finally, we prove Corollary \ref{coco} in
Section \ref{conseq}.

\section{Existence and uniqueness of the limit process}\label{ex}

The goal of this section is to show that the LFFP is well
defined,
unique and can be obtained from a graphical construction.
First, we show that when working on a finite
space interval, the LFPP is somewhat discrete.

We consider a Poisson measure $M(dt,dx)$ on $[0,\infty) \times
{\mathbb{R}}$ with
intensity measure $dt\, dx$. We define
${\mathcal F}_t^{M,A}=\sigma(M(B),B\in{\mathcal B}([0,t] \times[-A,A]))$.
\begin{defin}\label{dflffpA}
A $({\mathcal F}_t^{M,A})_{t\geq0}$-adapted process
\[
(Z_t^A(x),D_t^A(x), H_t^A(x))_{t\geq0,x\in[-A,A]}
\]
with values in
${\mathbb{R}}_+\times{\mathcal I}\times{\mathbb{R}}_+$
is called an $A$-LFFP
if a.s., for all $t\geq0$ and all $x \in[-A,A]$,
\[
\cases{
\displaystyle Z_t^A(x)= {\int_0^t}{\mathbf{1}}_{\{Z_s^A (x) < 1\}}\,ds -
{\int_0^t}\int_{[-A,A]} {\mathbf{1}}_{\{ Z_{{s-}}^A(x)=1,y \in
D_{{{s-}}}^A(x)\}}M(ds,dy),\cr
\displaystyle H_t^A(x)=
{\int_0^t}Z_{{s-}}^A(x){\mathbf{1}}_{\{Z_{{s-}}^A(x)<1\}} M(ds\times
\{x\})
- {\int_0^t}{\mathbf{1}}_{\{H_s^A (x) > 0 \}}\,ds,}
\]
where $D_t^A(x) = [L_t^A(x),R_t^A(x)]$ with
%
%
\begin{equation}\label{tructruc}
\cases{\displaystyle L_t^A(x) = (-A) \lor\sup\{ y \in[-A,x]; Z_t^A(y)<1
\mbox{ or } H_t^A(y)>0 \},\cr
\displaystyle R_t^A(x) = A \land\inf\{ y\in[x,A]; Z_t^A(y)<1 \mbox{ or } H_t^A(y)>0
\}.}
\end{equation}
\end{defin}

\begin{figure}

\includegraphics{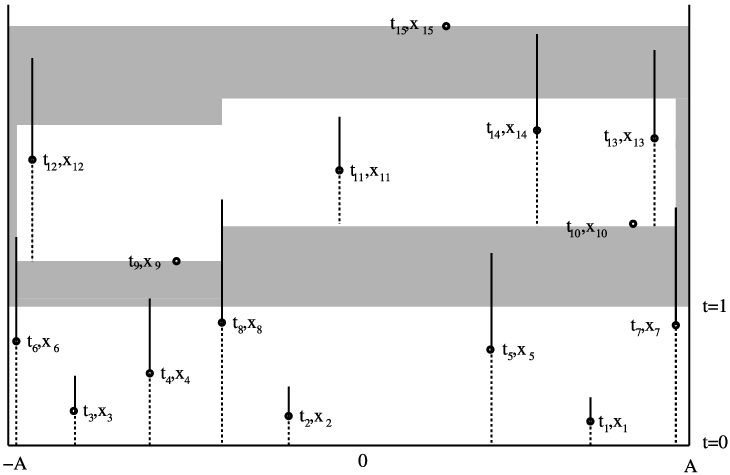}

\caption{Limit forest fire process in a finite box.
The filled zones represent zones in which $Z_t^A(x)=1$ and $H_t^A(x)=0$,
that is, macroscopic clusters. The plain vertical segments represent
the sites where $H_t^A(x)>0$. In the rest of the space, we always have
$Z_t^A(x)<1$.
Until time $1$, all of the particles are microscopic. The first eight
marks of the Poisson measure fall in that zone. As a consequence,
at each of these marks, the process $H^A$ starts. Their lifetime
is equal to the instant where they have started (e.g.,
the segment above $t_1,x_1$ ends at time $2t_1$).
At time $1$, all of the clusters where
there has been no mark become macroscopic and merge together.
However, this is limited by vertical segments. Here, at time $1$, we
have the
clusters $[-A,x_6]$, $[x_6,x_4]$, $[x_4,x_8]$, $[x_8,x_5]$, $[x_5,x_7]$ and
$[x_7,A]$. The segment above $(t_4,x_4)$ ends at time $2t_4$ and thus,
at this
time, the clusters $[x_6,x_4]$ and $[x_4,x_8]$ merge into
$[x_6,x_8]$.
The ninth mark
falls in the (macroscopic) zone $[x_6,x_8]$ and thus destroys it immediately.
This zone $[x_6,x_8]$ will become macroscopic again only at time
$t_9+1$. A process $H^A$ then starts at $x_{12}$ at time $t_{12}$.
Since\vspace*{1pt}
$Z^A_{t_{12}-}(x_{12})=t_{12}-t_9$
[because $Z^A_{t_9}(x_{12})$ has been set to $0$],
the segment above $(t_{12},x_{12})$ will end at time $2t_{12}-t_9$.
On the other hand, the segment $[x_8,x_7]$ has been destroyed at time $t_{10}$
and will thus remain microscopic until $t_{10}+1$. As a consequence,
the only macroscopic clusters at time $t_9+1$ are $[-A,x_{12}]$,
$[x_{12},x_8]$ and $[x_7,A]$. The zone $[x_8,x_7]$ then becomes macroscopic
(but there have been marks at $x_{13}, x_{14}$) so that at time $t_{10}+1$,
we get the macroscopic clusters $[-A,x_{12}]$, $[x_{12},x_{14}]$,
$[x_{14},x_{13}]$ and $[x_{13},A]$. These clusters merge by pairs, at times
$2t_{12}-t_9$,
$2t_{13}-t_{10}$ and $2t_{14}-t_{10}$,
so that we have a unique cluster $[-A,A]$
just before time $t_{15}$, where a mark falls and destroys the whole cluster
$[-A,A]$.\break
\mbox{\quad}With this realization, we have $0 \in(x_{11},x_{15})$ and, thus,
$Z^A_t(0)= t$ for $t \in[0,1]$, then $Z^A_t(0)=1$ for $t \in
[1,t_{10})$, then
$Z^A_t(0)= t-t_{10}$ for $t \in[t_{10},t_{10}+1)$, then
$Z^A_t(0)=1$ for $t \in[t_{10}+1,t_{15})$, etc.
We also see that $D^A_t(0)=\{0\}$ for $t \in[0,1)$,
$D^A_t(0)=[x_8,x_5]$ for $t\in[1,2t_5)$,
$D^A_t(0)=[x_8,x_7]$ for $t\in[2t_5,t_{10})$,
$D^A_t(0)=\{0\}$ for $t\in[t_{10},t_{10}+1)$, $D^A_t(0)=[x_{12},x_{14}]$
for $t \in[t_{10}+1,2t_{12}-t_9)$, $D^A_t(0)=[-A,x_{14}]$
for $t \in[2t_{12}-t_9,2t_{14}-t_{10})$, etc.
Of course, $H^A_t(0)=0$ for all $t\geq0$, but, for example,
$H^A_t(x_{11})=0$ for
$t\in[0, t_{11})$, $H^A_t(x_{11})=2t_{11}-t_{10}-t$ for
$t\in[t_{11},2t_{11}-t_{10})$
and then $H^A_t(x_{11})=0$ for $t\in[2t_{11}-t_{10},\infty)$.}
\label{figure2}
\end{figure}

A typical path of $(Z_t^A(x),D_t^A(x),H_t^A(x))_{t\geq0,x\in[-A,A]}$
is drawn in Figure \ref{figure2}.

%

Although the following proposition is almost obvious, its proof
shows the construction of the $A$-LFFP in an algorithmic way.
\begin{prop}\label{eufini}
Consider a Poisson measure $M(dt,dx)$ on $[0,\infty) \times{\mathbb
{R}}$ with
intensity measure $dt\, dx$. For any $A>0$, there a.s. exists a unique
$A$-LFFP which can be perfectly simulated.
\end{prop}
\begin{pf} We omit the superscript $A$ in this proof.
We consider the marks $(T_i,X_i)_{i \geq1}$ of $M\vert_{[0,\infty
)\times[-A,A]}$,
where $0<T_1<T_2<\cdots.$ We set $T_0=0$ for convenience.
We describe the construction via an algorithm,
which also shows uniqueness, in the sense that there is no choice in
the construction.\vspace*{8pt}

\textit{Step} 0. First, we set $Z_0(x)=H_0(x)=0$ and $D_0(x)=\{x\}$
for all $x\in[-A,A]$.\vspace*{8pt}

\textit{Step} $n+1$. Assume that the process has been built until $T_n$
for some $n\geq0$, that is, we know the values of
$(Z_t(x),D_t(x),H_t(x))_{t\in[0,T_n],x\in[-A,A]}$.

We build $(Z_t(x),D_t(x),H_t(x))_{t\in(T_n,T_{n+1}),x\in[-A,A]}$
in the following way: for $t\in(T_n,T_{n+1})$ and $x\in[-A,A]$,
we set $Z_t(x)= \min(1,Z_{T_n}(x)+t-T_n)$, $H_t(x)=\max
(0,H_{T_n}(x)-(t-T_n))$ and
define $D_t(x)=[L_t(x),R_t(x)]$, as in (\ref{tructruc}).

Next, we build
$(Z_{T_{n+1}}(x),D_{T_{n+1}}(x),H_{T_{n+1}}(x))_{x\in[-A,A]}$.

{\smallskipamount=0pt
\begin{longlist}
\item
If $Z_{T_{n+1}-}(X_{n+1})=1$, then we set $H_{T_{n+1}}(x)
=H_{T_{n+1}-}(x)$ for all
$x\in[-A,A]$ and consider $[a,b]:= D_{T_{n+1}-}(X_{n+1})$.
Set $Z_{T_{n+1}}(x) = 0$ for all $x\in(a,b)$ and $Z_{T_{n+1}}(x) =
Z_{T_{n+1}-}(x)$
for all $x\in[-A,A] \setminus[a,b]$. Finally, set:\break
$Z_{T_{n+1}}(a)=0$ if $Z_{T_{n+1}-}(a)=1$; $Z_{T_{n+1}}(a)=Z_{T_{n+1}-}(a)$
if $Z_{T_{n+1}-}(a)<1$;\break
$Z_{T_{n+1}}(b)=0$ if $Z_{T_{n+1}-}(b)=1$; $Z_{T_{n+1}}(b)=Z_{T_{n+1}-}(b)$
if $Z_{T_{n+1}-}(b)<1$.

\item If $Z_{T_{n+1}-}(X_{n+1})<1$, then we set
$H_{T_{n+1}}(X_{n+1})=Z_{T_{n+1}-}(X_{n+1})$,\break
$Z_{T_{n+1}}(X_{n+1})=Z_{T_{n+1}-}(X_{n+1})$ and
$(Z_{T_{n+1}}(x),H_{T_{n+1}}(x))
=(Z_{T_{n+1}-}(x)$,\break $H_{T_{n+1}-}(x))$ for
all $x\in[-A,A]\setminus\{X_{n+1}\}$.

\item Using the values of $(Z_{T_{n+1}}(x),H_{T_{n+1}}(x))_{x\in[-A,A]}$,
we finally compute the values of $(D_{T_{n+1}}(x))_{x\in[-A,A]}$.\qed
\end{longlist}}
\noqed\end{pf}

In case (i) above, we explained precisely what is done at the boundary of
burning macroscopic components. This is not so important: it does not
affect the uniqueness statement, but corresponds to using a slightly
different definition of the process; we could have made
other choices for this.

We now prove a refined version of Theorem \ref{existe}.
\begin{prop}\label{loc}
Consider a Poisson measure $M(dt,dx)$ on $[0,\infty) \times{\mathbb
{R}}$ with
intensity measure $dt \,dx$. For $A>0$, consider the
$A$-LFFP $(Z^A_t(x),D^A_t(x)$,\break $H^A_t(x))_{t\geq0,x\in[-A,A]}$
constructed in Proposition \ref{eufini} (using $M$).

There a.s. exists a unique LFFP
$(Z_t(x),D_t(x),H_t(x))_{t\geq0,x\in{\mathbb{R}}}$ (corresponding to $M$)
and, furthermore, it is such that for all $T>0$, there are constants
$\alpha_T>0$ and $C_T>0$ such that for all $A\geq2$,
%
%
\begin{eqnarray}\label{mix}
&&{\mathbb{P}}\bigl[(Z_t(x),D_t(x),H_t(x))_{t\in[0,T],x\in[-A/2,A/2]}
\nonumber\\[-8pt]\\[-8pt]
&&\qquad= (Z_t^A(x),D_t^A(x),H_t^A(x))_{t\in[0,T],x\in[-A/2,A/2]} \bigr]
\geq1 - C_T e^{-\alpha_T A}.\nonumber
\end{eqnarray}
\end{prop}
\begin{pf} We divide the proof into several steps. We fix $T>0$
and work on $[0,T]$.\vspace*{8pt}

\textit{Step} 1. For $a \in\mathbb{Z}$, we define the event
$\Omega_a$ in the following way (see Figure \ref{figure3} for an
illustration).
%
%
\begin{figure}[b]

\includegraphics{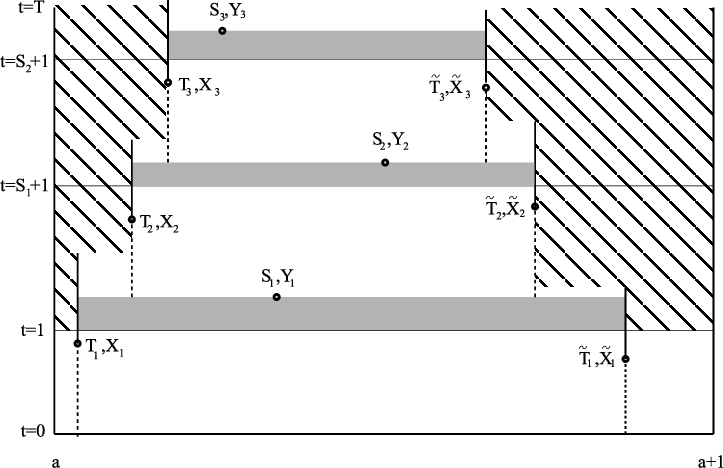}

\caption{The event $\Omega_a$ (proof of Theorem \protect\ref{existe}).
In hatched zones, we cannot state the values of the LFFP
because one would need to know what happens outside $[a,a+1]$.
\break
\mbox{\quad}Microscopic fires start at $(T_1,X_1)$ and $({\tilde T}_1,{\tilde X}_1)$.
Hence, at time $S_1$, the connected component $[X_1,{\tilde X}_1]$
is macroscopic because $S_1\geq1$ and because during $[1,S_1)$, this
component
has not been subject to fires starting outside $[a,a+1]$: it is
protected by
$X_1$ and ${\tilde X}_1$ until time $2\min(T_1,{\tilde T}_1) \geq S_1$.
As a consequence, the component $[X_1,{\tilde X}_1]$ is entirely killed by
$(S_1,Y_1)$.
We then iterate the arguments until we reach the final time $T$.
\break
\mbox{\quad}With such a configuration, there are always \textit{microscopic} sites in $[a,a+1]$
during $[0,T]$. Indeed, during $[0,1)$, all of the sites are microscopic,
during $[1,S_1)$,
the sites $X_1$ and ${\tilde X}_1$ are microscopic, during
$[S_1,S_1+1)$, all
the sites in $[X_1,{\tilde X}_1]$
are microscopic, etc.}
\label{figure3}
\end{figure}
The Poisson measure $M$
has exactly $3n$ marks in $[0,T]\times[a,a+1]$ for some
$n\geq1$ and
it is possible to call them $(T_k,X_k)_{k=1,\ldots,n}$,
$({\tilde T}_k,{\tilde X}_k)_{k=1,\ldots,n}$
and $(S_k,Y_k)_{k=1,\ldots,n}$
in such a way that
we have the following properties for all $k=1,\ldots,n$
(we set $T_0={\tilde T}_0=S_0=0$ and $X_0=a$,
${\tilde X}_0=a+1$ for convenience):

{\smallskipamount=0pt
\begin{longlist}
\item
$T_k$ and ${\tilde T}_k$ belong to $(S_{k-1}+1/2,S_{k-1}+1)$
and $X_{k-1}<X_k<{\tilde X}_k<{\tilde X}_{k-1}$;

\item $S_k\in(S_{k-1}+1,S_{k-1}+2(T_k\land{\tilde T}_k-S_{k-1}))$ and
$Y_k \in(X_k,{\tilde X}_k)$;

\item $S_n > T-1$.
\end{longlist}}
\vspace*{8pt}

\textit{Step} 2. We next observe that if the LFFP exists, then,
necessarily,
\[
\Omega_a \subset\{\forall t\in[0,T], \exists x\in(a,a+1),
H_t(x)>0 \mbox{ or } Z_t(x)<1 \}.
\]
Indeed, $Z_t(x) =t <1$ for all $t\in[0,1)$ and $x\in{\mathbb{R}}$.
Then $H_{T_1}(X_1)=Z_{T_1}(X_1)=T_1$, whence\vspace*{1pt} $H_t(X_1)>0$ on $[T_1,
2T_1]$ and
$H_t({\tilde X}_1)>0$ on $[{\tilde T}_1,2{\tilde T}_1]$. As a
consequence, we know\vspace*{1pt}
that for all $x\in(X_1,{\tilde X}_1)$ and $t\in[1,S_1)$, we have
$D_t(x)=[X_1,{\tilde X}_1]$. Since, now, $1<S_1<2(T_1\land{\tilde
T}_1)$ and since
$Y_1 \in(X_1,{\tilde X}_1)$, we deduce that $Z_{S_1}(x)=0$ for all
$x\in(X_1,{\tilde X}_1)$ and, as a consequence, $Z_t(x)=t-S_1<1$ for all
$t\in[S_1,S_1+1)$. However, we now have $H_t(X_2)>0$ on $[T_2,
T_2+(T_2-S_1))$ and
$H_t({\tilde X}_2)>0$ on $[{\tilde T}_2,{\tilde T}_2+({\tilde
T}_2-S_1))$. As a consequence, we know
for all $x\in(X_2,{\tilde X}_2)$ and $t\in[S_1+1,S_2)$ that
$D_t(x)=[X_2,{\tilde X}_2]$. Since, now, $S_1+1<S_2<S_1+2(T_1\land
{\tilde T}_1-S_1)$
and $Y_2 \in(X_2,{\tilde X}_2)$, we deduce that $Z_{S_2}(x)=0$ for all
$x\in(X_2,{\tilde X}_2)$ and thus $Z_t(x)=t-S_2<1$ for all
$t\in[S_2,S_2+1)$, etc.\vspace*{8pt}

\textit{Step} 3. We deduce that for all $a\in\mathbb{Z}$,
conditionally on $\Omega_a$,
clusters to the left of $a$ are never connected (during $[0,T]$)
to clusters to the right of $a+1$. Thus, on $\Omega_a$,
fires starting to the left of
$a$ do not affect the zone $[a+1,\infty)$ and fires starting to
the right of $a+1$ do not affect the zone $(-\infty,a]$.
Since, further, $\Omega_a$ concerns the Poisson measure $M$ only
in $[0,T]\times[a,a+1]$, we deduce that on $\Omega_a$,
the processes $(Z_t(x),D_t(x),H_t(x))_{t\geq0,x\in[a+1,\infty)}$
and $(Z_t(x),D_t(x),H_t(x))_{t\geq0,x\in(-\infty, a]}$
can be constructed separately.\vspace*{8pt}

\textit{Step} 4. Clearly, $q_T={\mathbb{P}}[\Omega_a]$ does not depend
on $a$,
by translation invariance (of the law of $M$), and obviously $q_T>0$.
Thus, a.s. there are infinitely many $a\in\mathbb{Z}$ such that
$\Omega_a$
is realized. This allows a graphical construction: it suffices to work
between such $a$'s (i.e., in finite boxes), as
in Proposition \ref{eufini}.\vspace*{8pt}

\textit{Step} 5. Using the same arguments, we easily deduce that
for $A\geq2$, the LFFP and the $A$-LFFP coincide on $[-A/2,A/2]$
during $[0,T]$, provided that there are $a_1 \in[-A,-A/2-1]$ and
$a_2\in[A/2,A-1]$
with $\Omega_{a_1} \cap\Omega_{a_2}$ realized.
Furthermore, since $M$ is a Poisson measure, $\Omega_a$ is independent
of $\Omega_b$ for all $a\ne b$ (with \mbox{$a,b\in\mathbb{Z}$}).
Thus, the
probability on the left-hand side of (\ref{mix}) is bounded below,
for $A\geq2$, by
\[
1- {\mathbb{P}}\biggl[\bigcap_{a\in\mathbb{Z}\cap[-A,-A/2-1]} \Omega_a^c\biggr] -
{\mathbb{P}}\biggl[\bigcap_{a\in\mathbb{Z}\cap[A/2, A-1]} \Omega_a^c\biggr] \geq
1 - 2 (1-q_T)^{A/2-2},
\]
hence we have (\ref{mix})
with $\alpha_T=-\log(1-q_T)/2>0$ and $C_T=2/(1-q_T)^2$.
\end{pf}

\section{Localization of the FFP}\label{locla}

We first introduce the $({\lambda},A)$-FFP.
We consider two independent families of i.i.d. Poisson processes
$N=(N_t(i))_{t\geq0,i\in\mathbb{Z}}$ and
$M^{\lambda}=(M^{\lambda}_t(i))_{t\geq0,i\in\mathbb{Z}}$, with
respective rates
$1$ and ${\lambda}>0$.
For $A>0$ and ${\lambda}>0$,
we define
%
%
\begin{equation}\label{alaila}
A_{\lambda}:=\bigl\lfloor A/\bigl({\lambda}\log(1/{\lambda})\bigr)\bigr\rfloor
\quad\mbox{and}\quad
I_A^{\lambda}:=\llbracket-A_{\lambda},A_{\lambda}\rrbracket,
\end{equation}
and we set
${\mathcal F}_t^{N,M^{\lambda},A}:=\sigma(N_s(i),M^{\lambda}_s(i),
s\leq t, i \in
I_A^{\lambda})$.
\begin{defin}\label{dflaffpA}
Consider an
$({\mathcal F}_t^{N,M^{\lambda},A})_{t\geq0}$-adapted
process $(\eta^{{\lambda},A}_t)_{t\geq0}$ with values in $\{0,1\}
^{I_A^{\lambda}}$,
such that $(\eta^{{\lambda},A}_t(i))_{t\geq0}$
is a.s. c\`adl\`ag for all $i\in I_A^{\lambda}$.

We say that $(\eta^{{\lambda},A}_t)_{t\geq0}$
is a $({\lambda},A)$-FFP if a.s., for all $t\geq0$ and $i\in
I_A^{\lambda}$,
\[
\eta^{{\lambda},A}_t(i)={\int_0^t}{\mathbf{1}}_{\{\eta^{{\lambda
},A}_{{s-}}(i)=0\}} \,dN_s(i)
- \sum_{k\in I_A^{\lambda}}{\int_0^t}{\mathbf{1}}_{\{k\in
C^{{\lambda},A}_{{s-}}(i)\}}
\,dM^{\lambda}_s(k),
\]
where $C^{{\lambda},A}_s(i)=\varnothing$ if $\eta_t^{{\lambda
},A}(i)=0$, while
$C^{{\lambda},A}_s(i)=\llbracket l_s^{{\lambda},A}(i),r_s^{{\lambda
},A}(i)\rrbracket$
if $\eta^{{\lambda},A}_s(i)=1$, where
\begin{eqnarray*}
l_s^{{\lambda},A}(i)&=&(-A_{\lambda}) \lor\bigl(\sup\{k< i; \eta
_s^{{\lambda},A}(k)=0\}+1
\bigr),\\
r_s^{{\lambda},A}(i)&=&A_{\lambda}\land\bigl(\inf\{k > i; \eta
_s^{{\lambda},A}(k)=0\}-1 \bigr).
\end{eqnarray*}
For $x \in[-A,A]$ and $t\geq0$, we introduce
%
%
\begin{eqnarray}
\label{dlambdaA}
D^{{\lambda},A}_t(x)&=&{\lambda}\log(1/{\lambda}) C^{{\lambda
},A}_t\bigl(\bigl\lfloor x/\bigl({\lambda}\log
(1/{\lambda})\bigr)\bigr\rfloor\bigr)
\subset[-A,A],
\\
%
%
\label{zlambdaA}
Z^{{\lambda},A}_t(x)&=&\frac{\log
[1+\# (C^{{\lambda},A}_t(\lfloor x/({\lambda}\log(1/{\lambda
}))\rfloor) ) ]}
{\log(1/{\lambda})} \geq0.
\end{eqnarray}
\end{defin}

We now prove the following result, which
is similar to Proposition \ref{loc} for the ${\lambda}$-FFP.
\begin{prop}\label{plocla}
Let $T>0$ and ${\lambda}\in(0,1)$. Consider two families of Poisson processes
$N=(N_t(i))_{t\geq0,i\in\mathbb{Z}}$ and
$M^{\lambda}=(M^{\lambda}_t(i))_{t\geq0,i\in\mathbb{Z}}$ with
respective rates
$1$ and ${\lambda}>0$. Let $(\eta^{\lambda}_t)_{t\geq0}$ be the
corresponding
${\lambda}$-FFP
and, for each $A>0$, let
$(\eta^{{\lambda},A}_t)_{t\geq0}$ be the corresponding $({\lambda},A)$-FFP.
Recall (\ref{dlambda}), (\ref{zlambda}) and (\ref{dlambdaA}), (\ref
{zlambdaA}).
There are constants $\alpha_T>0$ and $C_T>0$,
not depending on ${\lambda}\in(0,1)$, $A\geq2$, such that [recalling
(\ref
{alaila})]
\begin{eqnarray*}
&&{\mathbb{P}}\bigl[(\eta^{\lambda}_t(i))_{t\in[0,T\log(1/{\lambda
})],i\in I_{A/2}^{\lambda}}
=(\eta^{{\lambda},A}_t(i))_{t\in[0,T\log(1/{\lambda})],i\in
I_{A/2}^{\lambda}} \bigr]\\
&&\qquad
\geq1 - C_T e^{-\alpha_T A},\\
&&{\mathbb{P}}\bigl[(Z^{\lambda}_t(x),D^{\lambda}_t(x))_{t\in[0,T],x \in
[-A/2,A/2]}
=(Z^{{\lambda},A}_t(x),D^{{\lambda},A}_t(x))_{t\in[0,T],x \in
[-A/2,A/2]} \bigr] \\
&&\qquad \geq1 - C_T e^{-\alpha_T A}.
\end{eqnarray*}
\end{prop}
\begin{pf} The proof is similar (but more
complicated) to that of Proposition~\ref{loc}.
Consider the true ${\lambda}$-FFP $(\eta_t^{\lambda}(i))_{t\geq0, i
\in\mathbb{Z}}$.
Temporarily assume that for $a \in{\mathbb{R}}$,
there is an event $\Omega_a^{\lambda}$, depending only on
the Poisson processes $N_t(i)$ and $M^{\lambda}_t(i)$ for $t\in
[0,T\log
(1/{\lambda})]$
and $i \in J_a^{\lambda}:= \llbracket\lfloor a/({\lambda}\log
(1/{\lambda}))\rfloor,
\lfloor(a+1)/({\lambda}\log(1/{\lambda}))\rfloor\rrbracket$, such that:

{\smallskipamount=0pt
\begin{longlist}
\item
on $\Omega^{\lambda}_a$, a.s., for all $t\in[0,T\log(1/{\lambda
})]$, there
is some
$i \in J_a^{\lambda}$ such that $\eta^{\lambda}_t(i)=0$;

\item there exists $q_T>0$ such that for all $a \in{\mathbb{R}}$ and
${\lambda}\in
(0,1)$, we have
${\mathbb{P}}(\Omega^{\lambda}_a) \geq q_T$.
\end{longlist}}

The proof is then concluded using arguments similar to Steps 3, 4, 5
of the proof of Proposition \ref{loc}.

Fix some $\alpha>0$ and some ${\varepsilon}_T>0$ small enough, say
$\alpha=0.01$ and
${\varepsilon}_T = 1/(32T)$.
Let ${\lambda}_T>0$ be such that for ${\lambda}\in(0,{\lambda}_T)$,
we have
$1<{\lambda}^{\alpha-1} < \epsilon_T /\break({\lambda}\log(1/{\lambda}))$.

For ${\lambda}\in[{\lambda}_T,1)$ and $a\in{\mathbb{R}}$, we set
$\Omega_a^{\lambda}=\{N_{T\log(1/{\lambda})}(\lfloor a/({\lambda
}\log(1/{\lambda}
))\rfloor)=0\}$,
on which, of course,\vspace*{1pt}
$\eta_t^{\lambda}(i)=0$ for all $t\in[0,T\log(1/{\lambda})]$ with
$i=\lfloor a/({\lambda}\log(1/{\lambda}))\rfloor\in J_a^{\lambda}$.
We then observe that $q'_T= \inf_{{\lambda}\in[{\lambda}_T,1)}
P(\Omega_a^{\lambda})
= \inf_{{\lambda}\in[{\lambda}_T,1)} e^{-T\log(1/{\lambda})}=
({\lambda}_T)^T>0$.

For ${\lambda}\in(0,{\lambda}_T)$ and $a\in{\mathbb{R}}$, we
define the event
$\Omega^{\lambda}_a$ on which points 1, 2 and 3 below are
satisfied.

1. The family of Poisson processes $(M^{\lambda}_t(i))_{t \in[0,T
\log
(1/{\lambda})],
i \in J_a^{\lambda}}$ has exactly $3n$ marks for some
$1 \leq n \leq\lfloor T \rfloor$ and
it is possible to call them $(T^{\lambda}_k,X^{\lambda}_k)_{k=1,\ldots,n}$,
$({\tilde T}^{\lambda}_k,{\tilde X}^{\lambda}_k)_{k=1,\ldots,n}$
and $(S^{\lambda}_k,Y^{\lambda}_k)_{k=1,\ldots,n}$
in such a way that
we have the following properties for all
$k=1,\ldots,n$ (we set $T^{\lambda}_0={\tilde T}^{\lambda}_0=S^{\lambda
}_0=0$ and
$X^{\lambda}_0=\lfloor a/\break({\lambda}\log(1/{\lambda}))\rfloor$,
${\tilde X}^{\lambda}_0=\lfloor(a+1)/({\lambda}\log(1/{\lambda
}))\rfloor$):

(1a) $X^{\lambda}_{k-1} < X^{\lambda}_k < Y^{\lambda}_k< {\tilde
X}^{\lambda}_k< {\tilde X}^{\lambda}_{k-1}$ with
$\min\{ X^{\lambda}_k-X^{\lambda}_{k-1}, Y^{\lambda}_k - X^{\lambda
}_k, {\tilde X}^{\lambda}_k -
Y^{\lambda}_k,
{\tilde X}^{\lambda}_{k-1} - {\tilde X}^{\lambda}_k\} \geq4 \epsilon
_T / ({\lambda}\log(1/{\lambda}))$;

(1b) $T^{\lambda}_k$ and ${\tilde T}^{\lambda}_k$ belong to
$[S^{\lambda}_{k-1}+(\frac
{1}{2} +
\alpha) \log(1/{\lambda}), S^{\lambda}_{k-1}+(1-\alpha) \log
(1/{\lambda})]$;

(1c) $S^{\lambda}_k\in[S^{\lambda}_{k-1}+(1+ \alpha) \log
(1/{\lambda}), S^{\lambda}_{k-1}
+2(T^{\lambda}_k\land{\tilde T}^{\lambda}_k-S^{\lambda}_{k-1}) -
\alpha\log(1/{\lambda})]$;

(1d) $S^{\lambda}_n\geq(T-1+\alpha) \log(1/{\lambda})$.

2. For $k=1,\ldots,n$, we now set
$\tau_k^{\lambda}=(S^{\lambda}_k-S^{\lambda}_{k-1})/(2\log
(1/{\lambda}))$, which belongs
to $[(1+\alpha)/2,1-\alpha]$, due to 1. We consider the intervals
\begin{eqnarray*}
I^{\lambda}_k&=&\llbracket X^{\lambda}_k- \lfloor{\lambda}^{-\tau
_k^{\lambda}} \rfloor,
X^{\lambda}_k
+ \lfloor{\lambda}^{-\tau_k^{\lambda}} \rfloor\rrbracket
,\nonumber\\
I^{\lambda}_{k,-}&=& \llbracket X^{\lambda}_k - \lfloor{\lambda
}^{-\tau_k^{\lambda}}
\rfloor
- \lfloor{\varepsilon}_T/{\lambda}\log(1/{\lambda}) \rfloor,
X^{\lambda}_k -
\lfloor{\lambda}^{-\tau_k^{\lambda}} \rfloor-1 \rrbracket
,\nonumber\\
I^{\lambda}_{k,+}&=& \llbracket X^{\lambda}_k + \lfloor{\lambda
}^{-\tau_k^{\lambda}}
\rfloor+1 ,
X^{\lambda}_k + \lfloor{\lambda}^{-\tau_k^{\lambda}} \rfloor+
\lfloor{\varepsilon
}_T/{\lambda}\log
(1/{\lambda})
\rfloor\rrbracket,\nonumber\\
L^{\lambda}_k&=& \llbracket X^{\lambda}_k + \lfloor{\lambda}^{-\tau
_k^{\lambda}} \rfloor+
\lfloor
{\varepsilon}_T/{\lambda}\log(1/{\lambda}) \rfloor+1,\\
&&\hspace*{8.3pt}
{\tilde X}^{\lambda}_k - \lfloor{\lambda}^{-\tau_k^{\lambda}}
\rfloor-
\lfloor{\varepsilon}_T/{\lambda}\log(1/{\lambda}) \rfloor-1
\rrbracket
\end{eqnarray*}
and similar intervals ${\tilde I}^{\lambda}_k,{\tilde I}^{\lambda
}_{k,-},{\tilde I}^{\lambda}_{k,+}$,
around ${\tilde X}^{\lambda}_k$.
For all $k=1, \ldots, n$, the family of Poisson processes
$(N_t(i))_{t\geq0, i \in J^{\lambda}_a}$ satisfies:

(2a) $\forall i \in I_k^{\lambda}, N_{T^{\lambda}_{k}}(i) -
N_{S^{\lambda}_{k-1}}(i)
>0$ and
$\forall i \in{\tilde I}_k^{\lambda}, N_{{\tilde T}^{\lambda
}_{k}}(i) - N_{S^{\lambda}_{k-1}}(i) >0$;

(2b) $\exists i \in I_{k,-}^{\lambda}$ such that $ N_{T^{\lambda
}_{k}}(i) -
N_{S^{\lambda}_{k-1}}(i) =0 $,
$\exists i \in I_{k,+}^{\lambda}$ such that $N_{T^{\lambda}_{k}}(i) -
N_{S^{\lambda}
_{k-1}}(i) =0 $,
$\exists i \in{\tilde I}_{k,-}^{\lambda}$ such that $N_{{\tilde
T}^{\lambda}_{k}}(i) -
N_{S^{\lambda}
_{k-1}}(i) =0 $
and $\exists i \in{\tilde I}_{k,+}^{\lambda}$ such that $N_{{\tilde
T}^{\lambda}_{k}}(i) -
N_{S^{\lambda}_{k-1}}(i) =0$;

(2c) $\exists i \in I_{k}^{\lambda}$ such that $N_{S^{\lambda
}_{k}}(i) - N_{T^{\lambda}
_{k}}(i) =0$ and
$\exists i \in{\tilde I}_{k}^{\lambda}$ such that $N_{S^{\lambda
}_{k}}(i) - N_{{\tilde T}
^{\lambda}
_{k}}(i) =0$;

(2d) $\forall i \in L_{k}^{\lambda}, N_{S^{\lambda}_{k}}(i) -
N_{S^{\lambda}_{k-1}}(i) >0$.

3. We finally assume that $\exists i \in L_{n}^{\lambda}$ such that
$N_{T\log(1/{\lambda})}(i) - N_{S^{\lambda}_{n}}(i)=0$.

To show that on $\Omega^{\lambda}_a$, a.s., for all $t\in[0,T\log
(1/{\lambda})]$,
there is some
$i \in J_a^{\lambda}$ such that $\eta^{\lambda}_t(i)=0$, we proceed
recursively.
At time $0$, all sites are vacant. Fix $k \in\{ 1,\ldots, n\}$. Assume
that for
$t \leq S^{\lambda}_{k-1}$, there is some
$i \in J_a^{\lambda}$ such that $\eta^{\lambda}_t(i)=0$ and that at time
$S^{\lambda}_{k-1}$, all sites in the interval $L_{k-1}^{\lambda}$
are vacant.

Then,
for $S^{\lambda}_{k-1} \leq t < T^{\lambda}_k$ (resp., $ S^{\lambda
}_{k-1} \leq t < {\tilde T}
^{\lambda}_k$),
(2b) shows that there are vacant sites in both $I_{k,+}^{\lambda}$ and
$I_{k,-}^{\lambda}$
(resp., in both ${\tilde I}_{k,+}^{\lambda}$ and ${\tilde
I}_{k,-}^{\lambda}$). This, together
with~(2a), shows that at time $T^{\lambda}_k-$ (resp., ${\tilde
T}^{\lambda}_k-$),
all of the sites in the intervals $I_k^{\lambda}$ and ${\tilde
I}_k^{\lambda}$ are occupied
(no fire may burn those sites because they are protected by the vacant
sites in $I_{k,+}^{\lambda}, I_{k,-}^{\lambda}, {\tilde
I}_{k,+}^{\lambda}, {\tilde I}_{k,-}^{\lambda}$).
Hence, the interval $I_k^{\lambda}$ (resp., ${\tilde I}_k^{\lambda}$)
becomes completely
vacant at time
$T^{\lambda}_k$ (resp., ${\tilde T}^{\lambda}_k$). Between time
$T^{\lambda}_k$ (resp., ${\tilde T}
^{\lambda}_k$)
and time~$S^{\lambda}_k$, since $I_k^{\lambda}$ (resp., ${\tilde
I}_k^{\lambda}$) is completely
vacant at
time $T^{\lambda}_k$ (resp., ${\tilde T}^{\lambda}_k$), (2c)~shows
that there is a vacant
site in $I_k^{\lambda}$ (resp., ${\tilde I}_k^{\lambda}$).

At time $S^{\lambda}_k-$, the interval $L_{k}^{\lambda}$ is
completely occupied,
by virtue of (2d) and the fact that it cannot be burnt because it is protected
by vacant sites in $I_{k,+}^{\lambda}$ (resp., ${\tilde
I}_{k,-}^{\lambda}$)
between $S^{\lambda}_{k-1}$ and $T^{\lambda}_{k}$ (resp., ${\tilde
T}^{\lambda}_k$), and in
$I_k^{\lambda}$
(resp., ${\tilde I}_k^{\lambda}$) between $T^{\lambda}_k$ (resp.,
${\tilde T}^{\lambda}_k$) and
$S^{\lambda}_k$.
As a consequence, since $Y^{\lambda}_k \in L_k^{\lambda}$, the interval
$L_k^{\lambda}$ becomes completely vacant at time $S^{\lambda}_k-$.

All of this shows that on $\Omega^{\lambda}_a$, there are vacant
sites in
$J^{\lambda}_a$
for all $t \in[0,S_n^{\lambda}]$ and that $L_n^{\lambda}$ is completely
vacant at time $S^{\lambda}_n$.
Finally, 3 implies that there are vacant sites in
$L_n^{\lambda}\subset J^{\lambda}_a$ during $[S_n^{\lambda},T\log
(1/{\lambda})]$.

It remains to prove that there exists $q''_T>0$ such that for all
$a \in{\mathbb{R}}$ and ${\lambda}\in(0,{\lambda}_T)$, we have
${\mathbb{P}}(\Omega^{\lambda}_a) \geq q''_T$. We separately treat
the conditions
1 on $M^{\lambda}$ and 2 on $N$
(conditionally on $M^{\lambda}$) and use independence
of these two families of Poisson processes to complete the proof.

First, for ${\lambda}\in(0, {\lambda}_T)$, we observe that we can
construct $M^{\lambda}$ using a Poisson measure $M$ on $[0,\infty
)\times
{\mathbb{R}}$
with intensity $dt\,dx$ by setting, for all
$i \in\mathbb{Z}$,
\[
M_t^{\lambda}(i)=M \bigl([0,t/\log(1/{\lambda})]\times
\bigl[i{\lambda}\log(1/{\lambda}), (i+1){\lambda}\log(1/{\lambda})\bigr) \bigr).
\]
Hence [since ${\varepsilon}_T/({\lambda}\log(1/{\lambda})) > 1$],
the event on which
$M^{\lambda}$ satisfies 1 contains the event $\Omega'_a$ on
which $M$
has exactly $3n$ marks in $[0,T]\times[a,a+1]$, for
some $1\leq n \leq\lfloor T \rfloor$, which can be called
$(T_k,X_k)_{k=1,\ldots,n}$, $({\tilde T}_k,{\tilde X}_k)_{k=1,\ldots,n}$
and $(S_k,Y_k)_{k=1,\ldots,n}$
in such a way that
we have the following properties (we set $T_0={\tilde T}_0=S_0=0$ and
$X_0=a$, ${\tilde X}_0=a+1$ for convenience) for all $k=1,\ldots,n$:

$\bullet$ $\min(\{ X_k-X_{k-1}, Y_k - X_k, {\tilde X}_k - Y_k,
{\tilde X}_{k-1} -
{\tilde X}
_k\})
> 5 \epsilon_T $;

$\bullet$ $T_k$ and ${\tilde T}_k$ belong to $(S_{k-1}+ 1/2 + \alpha,
S_{k-1}+ 1-\alpha);$

$\bullet$ $S_k \in(S_{k-1}+ 1+ \alpha, S_{k-1}+2(T_k \land{\tilde
T}_k-S_{k-1})
- \alpha) $;

$\bullet$ $S_n\geq(T-1) + \alpha$.

We then have ${\mathbb{P}}(\Omega_a')>0$ (as in the proof
of Proposition \ref{loc} and since ${\varepsilon}_T$ and $\alpha$
are sufficiently small)
and this probability does not depend on $a$ (by translation invariance
of the law of $M$) nor on ${\lambda}\in(0,{\lambda}_T)$ (since it
concerns only $M$).

We then use basic computations on i.i.d. Poisson processes
with rate $1$ to show that
there is a (deterministic) constant $c>0$ such that for all $k= 1,
\ldots, n$,
all ${\lambda}\in(0,\lambda_T)$, conditionally on $M^{\lambda}$
(we write ${\mathbb{P}}_M$ for the conditional probability w.r.t.
$M^{\lambda}$):

$\bullet$ since $T^{\lambda}_k - S^{\lambda}_{k-1}\geq(\tau
_k^{\lambda}+\alpha
/2)\log
(1/{\lambda})$,
due to (1c), and since $\#(I_k^{\lambda})=2\lfloor{\lambda}^{-\tau
_k^{\lambda}
}\rfloor+1$,
we have
\begin{eqnarray*}
{\mathbb{P}}_M\bigl( \forall i \in I_k^{\lambda}, N_{T^{\lambda}_{k}}(i)
- N_{S^{\lambda}
_{k-1}}(i) >0\bigr)
&=& \bigl(1-e^{-(T^{\lambda}_{k}-S^{\lambda}_{k-1})} \bigr)^{ 2\lfloor{\lambda
}^{-\tau_k^{\lambda}
}\rfloor
+1 }\\
&\geq&(1- {\lambda}^{\tau_k^{\lambda}+\alpha/2} )^{2\lfloor
{\lambda}^{-\tau_k^{\lambda}
}\rfloor+1}
\geq c
\end{eqnarray*}
(it tends to $1$ as ${\lambda}\to0$)
and the same computation works for ${\tilde I}^{\lambda}_k$;

$\bullet$ since $T^{\lambda}_k-S^{\lambda}_{k-1} \leq(1-\alpha)
\log(1/{\lambda}
)$, by (1b),
and since $\#(I_{k,+}^{\lambda})=\lfloor{\varepsilon}_T/({\lambda
}\times\break\log(1/{\lambda}
)) \rfloor$, we have
\begin{eqnarray*}
{\mathbb{P}}_M\bigl( \exists i \in I_{k,+}^{\lambda}, N_{T^{\lambda
}_{k}}(i) - N_{S^{\lambda}
_{k-1}}(i) =0\bigr)
&=& 1 - \bigl(1-e^{-(T^{\lambda}_k-S^{\lambda}_{k-1})}
\bigr)^{\lfloor{\varepsilon}_T /{\lambda}\log(1/{\lambda})\rfloor}
\nonumber\\
&\geq& 1- (1-{\lambda}^{1-\alpha} )^{\lfloor{\varepsilon
}_T/({\lambda}\log
(1/{\lambda}))\rfloor} \geq c
\end{eqnarray*}
and the same computation works for $I_{k,-}^{\lambda},{\tilde
I}_{k,+}^{\lambda},{\tilde I}
_{k,-}^{\lambda}$;

$\bullet$
since $S^{\lambda}_k - T^{\lambda}_k \leq(\tau_k^{\lambda}- \alpha
/2)\log(1/{\lambda}
)$, due
to (1c)
[we use the fact that $S_k^{\lambda}\leq2T_k^{\lambda
}-S_{k-1}^{\lambda}-\alpha\log
(1/{\lambda})$,
whence $2S_k^{\lambda}\leq2T_k^{\lambda}+S_k^{\lambda
}-S_{k-1}^{\lambda}-\alpha\log
(1/{\lambda})
=2T_k^{\lambda}+2 (\tau_k^{\lambda}-\alpha/2) \log(1/{\lambda})$],
and since $\#(I_k^{\lambda})=2\lfloor{\lambda}^{-\tau_k^{\lambda
}}\rfloor+1$,
we have
\begin{eqnarray*}
{\mathbb{P}}_M\bigl( \exists i \in I_k^{\lambda}, N_{S^{\lambda
}_{k}}(i) - N_{T^{\lambda}_{k}}(i)
=0\bigr) &=&
1 - \bigl(1-e^{- (S^{\lambda}_k-T^{\lambda}_k)} \bigr)^{2\lfloor{\lambda
}^{-\tau_k^{\lambda}
}\rfloor
+1}\\
&\geq& 1- (1-{\lambda}^{\tau_k^{\lambda}-\alpha/2} )^{2\lfloor
{\lambda}^{-\tau
_k^{\lambda}
}\rfloor+1}
\geq c
\end{eqnarray*}
and this also holds for ${\tilde I}^{\lambda}_k$;

$\bullet$ since $S^{\lambda}_k - S^{\lambda}_{k-1} \geq(1+\alpha)
\log(1/{\lambda})$,
thanks to (1c),
and since $\#(L_k^{\lambda}) \leq\lfloor(1/{\lambda}\log
(1/{\lambda}))\rfloor$, we have
\begin{eqnarray*}
{\mathbb{P}}_M\bigl( \forall i \in L_{k}^{\lambda}, N_{S^{\lambda
}_{k}}(i) - N_{S^{\lambda}
_{k-1}}(i) >0\bigr)
&=& \bigl( 1- e^{ -(S^{\lambda}_k - S^{\lambda}_{k-1}) } \bigr)^{\#(L_k^{\lambda})}
\\
&\geq& ( 1- {\lambda}^{1+\alpha} )^{\lfloor1/{\lambda}\log
(1/{\lambda})\rfloor}
\geq c;
\end{eqnarray*}

$\bullet$ since $T\log(1/{\lambda})-S_n^{\lambda}\leq(1-\alpha
)\log(1/{\lambda})$,
by (1d), and $\# (L_n^{\lambda}) \geq4{\varepsilon}_T/({\lambda
}\log(1/\break{\lambda}
))$, by (1a),
we have
\begin{eqnarray*}
{\mathbb{P}}_M\bigl(\exists i \in L_{n}^{\lambda},
N_{T\log(1/{\lambda})}(i) - N_{S^{\lambda}_{n}}(i)=0\bigr)&=&
1- \bigl(1-e^{-(T\log(1/{\lambda})-S_n^{\lambda})} \bigr)^{\# (L_n^{\lambda
})}\nonumber\\
&\geq& 1- (1-{\lambda}^{1-\alpha} )^{4{\varepsilon}_T/({\lambda}\log
(1/{\lambda}))}
\geq c.
\end{eqnarray*}

We observe that the domains $I_k^{\lambda}\times(S^{\lambda}_{k-1},
T^{\lambda}_k]$,
${\tilde I}^{\lambda}_k\times(S^{\lambda}_{k-1}, {\tilde T}^{\lambda}_k]$,
$I_{k,+}^{\lambda}\times(S^{\lambda}_{k-1}$,\break $T^{\lambda}_k]$,
$I_{k,-}^{\lambda}\times(S^{\lambda}_{k-1}, T^{\lambda}_k]$,
${\tilde I}_{k,+}^{\lambda}\times(S^{\lambda}_{k-1}, {\tilde
T}^{\lambda}_k]$,
${\tilde I}_{k,-}^{\lambda}\times(S^{\lambda}_{k-1}, {\tilde
T}^{\lambda}_k]$,
$I_k^{\lambda}\times(T^{\lambda}_{k}, S^{\lambda}_k]$,
${\tilde I}_k^{\lambda}\times({\tilde T}^{\lambda}_{k}, S^{\lambda}_k]$,
$L_{k}^{\lambda}\times(S^{\lambda}_{k-1}, S^{\lambda}_k]$, for
$k=1,\ldots, n$,
and $L_{n}^{\lambda}\times(S^{\lambda}_{n}, T\log(1/{\lambda})]$
are pairwise disjoint, thanks to 1 and to the smallness
of ${\varepsilon}_T$ and ${\lambda}_T$: we have $\lfloor{\lambda
}^{-\tau_k^{\lambda}
} \rfloor
\leq{\lambda}^{\alpha-1} \leq{\varepsilon}_T/({\lambda}\log
(1/{\lambda}))$.

Since $n\leq T$, we deduce from all of the previous
estimates the existence of a $q''_T>0$ such that for all
$a \in{\mathbb{R}}$ and ${\lambda}\in(0,{\lambda}_T)$, we have
${\mathbb{P}}(\Omega^{\lambda}_a) \geq q''_T$.
We complete the proof by choosing $q_T = \min(q'_T,q''_T)$.
\end{pf}

\section{Convergence proof}\label{cv}

The goal of this section is to prove Theorem \ref{converge}.

\subsection{Coupling}
We introduce a coupling between
the ${\lambda}$-FFP, the LFFP and their localized versions.
\begin{nota}\label{cpp}
We consider a Poisson measure $M(dt,dx)$ on $[0,\infty) \times
{\mathbb{R}}$
with intensity measure $dt\,dx$. We consider an independent family
of Poisson processes $(N_t(i))_{t\geq0, i\in\mathbb{Z}}$ with rate $1$.
For ${\lambda}\in(0,1)$ and $i \in\mathbb{Z}$, we set
\[
M_t^{\lambda}(i)=M \bigl([0,t/\log(1/{\lambda})]\times
\bigl[i{\lambda}\log(1/{\lambda}), (i+1){\lambda}\log(1/{\lambda})\bigr) \bigr).
\]
Then $(M_t^{\lambda}(i))_{t\geq0, i\in\mathbb{Z}}$ is a family of
independent
Poisson processes with rate ${\lambda}$.
For all ${\lambda}\in(0,1)$, we consider the ${\lambda}$-FFP
$(\eta^{\lambda}_t)_{t\geq0}$ (see Definition \ref{dflaffp}) and
for all $A>0$,
we consider the $({\lambda},A)$-FFP $(\eta^{{\lambda},A}_t)_{t\geq
0}$ (see
Definition \ref{dflaffpA})
constructed with $N,M^{\lambda}$. We also introduce the processes
$(Z^{\lambda}_t(x),D^{{\lambda}}_t(x))_{t\geq0,x\in{\mathbb{R}}}$,
as in (\ref{dlambda}), (\ref{zlambda}), and
$(Z^{{\lambda},A}_t(x),D^{{\lambda},A}_t(x))_{t \geq0,x\in[-A,A]}$,
as in (\ref{dlambdaA}), (\ref{zlambdaA}).

We denote by $(Z_t(x),D_t(x),H_t(x))_{t\geq0, x\in{\mathbb{R}}}$ the LFFP
constructed with
$M$ (see Definition \ref{dflffp}) and by
$(Z_t^A(x),D_t^A(x),H_t^A(x))_{t\geq0, x\in[-A,A]}$ the
$A$-LFFP constructed with $M$ (see Definition \ref{dflffpA}).
\end{nota}

\subsection{Localization}

Temporarily assume that the following result holds.
\begin{prop}\label{heart}
Adopt Notation \ref{cpp} as well as Notation \ref{nocv}.

\textup{(a)} For any $T>0$, $A>0$ and $x_0\in(-A,A)$, in probability,
as ${\lambda}\to0$,
\[
\delta_T ( (Z^{{\lambda},A}(x_0),D^{{\lambda},A}(x_0)),
(Z^{A}(x_0),D^A(x_0)) )\qquad
\mbox{tends to } 0.
\]

\textup{(b)} For any $t\in[0,\infty)$,
$A>0$ and $x_0\in(-A,A)$, in probability,
as ${\lambda}\to0$,
\[
|Z^{{\lambda},A}_t(x_0)-Z^{A}_t(x_0)|+\delta(D^{{\lambda
},A}_t(x_0)),D^A_t(x_0) )\qquad
\mbox{tends to } 0.
\]
\end{prop}

We are now in a position to give the
following proof.
\begin{pf*}{Proof of Theorem \ref{converge}}
We only prove point (a), (b) being similarly checked.
Let $T>0$ and $\{x_1,\ldots,x_n\} \subset[-B,B]\subset
{\mathbb{R}}$ be fixed.
Consider the coupling introduced in
Notation \ref{cpp}.
Proposition \ref{heart} ensures us that for any ${\varepsilon}>0$ and
$A>B$, we have
\[
\lim_{{\lambda}\to0}{\mathbb{P}}\Biggl[\sum_{1}^n \delta_T (
(Z^{{\lambda},A}(x_i),D^{{\lambda},A}(x_i)),
(Z^{A}(x_i),D^A(x_i)) ) >{\varepsilon}\Biggr] =0.
\]
Now, let
\begin{eqnarray*}
&&\Omega_{A,T}^{\lambda}:= \{\forall i =1,\ldots,n, \forall
t\in[0,T],\\
&&\hspace*{44.27pt}
(Z^{{\lambda}}_t(x_i),D^{{\lambda}}_t(x_i))=
(Z^{{\lambda},A}_t(x_i),D^{{\lambda},A}_t(x_i)) \\
&&\hspace*{44.27pt} \mbox{and }(Z_t(x_i),D_t(x_i))= (Z_t^A(x_i),D_t^A(x_i)) \}.
\end{eqnarray*}
For all $A>2B$,
we now have
\begin{eqnarray*}
&&\Omega_{A,T}^{\lambda}\subset\bigl\{
(Z_t^{\lambda}(x),D_t^{\lambda}(x))_{t\in[0,T],x\in[-A/2,A/2]}
\\
&&\hspace*{42.5pt}=
(Z_t^{{\lambda},A}(x),D_t^{{\lambda},A}(x))_{t\in[0,T],x\in
[-A/2,A/2]} \\
&&\hspace*{42.5pt}\mbox{and } (Z_t(x),D_t(x))_{t\in[0,T],x\in[-A/2,A/2]}\\
&&\hspace*{42.5pt}\hspace*{14.29pt}=
(Z_t^A(x),D_t^A(x))_{t\in[0,T],x\in[-A/2,A/2]} \bigr\}.
\end{eqnarray*}
However, Propositions \ref{loc} and \ref{plocla} yield that
${\mathbb{P}}[(\Omega_{A,T}^{\lambda})^c] \leq2C_T e^{-\alpha_T
A}$. Thus, for any
$A>2B$,
\[
\limsup_{{\lambda}\to0}{\mathbb{P}}\Biggl[\sum_{1}^n \delta_T
((Z^{{\lambda}}(x_i),D^{{\lambda}}(x_i)),
(D(x_i),Z(x_i)) ) >{\varepsilon}\Biggr] \leq0 + 2C_T e^{-\alpha_T A}.
\]
Letting $A$ tend to infinity, we deduce that
$\sum_{i=1}^n\delta_T((Z^{{\lambda}}(x_i),D^{{\lambda}}(x_i)),
(D(x_i)$,\break $Z(x_i)))$ tends to $0$ in probability as ${\lambda}\to0$,
hence the result.
\end{pf*}

\subsection{Core of the proof}

The aim of this subsection is to prove Proposition
\ref{heart}.
We fix $T>0$ and $A>0$. We consider the
$({\lambda},A)$-FFP and the $A$-LFFP coupled, as in Notation \ref{cpp},
and use the notation introduced in (\ref{alaila}).
Throughout this proof, we will omit the superscript $A$ and we do not
take into account the possible dependencies in $A$ and $T$.

For $J=(a,b)$ [an open interval of $(-A,A)$], ${\lambda}\in(0,1)$ and
$\mu
\in(0,1]$, we consider
%
%
\begin{eqnarray}\label{xxbb}
J_{{\lambda},\mu} &=& \biggl[\hspace*{-3.5pt}\biggl[\biggl\lfloor\frac{a}{{\lambda}\log
(1/{\lambda})}+\frac{\mu
}{{\lambda}\log
^2(1/{\lambda})}
\biggr\rfloor, \nonumber\\
&&\hspace*{9.9pt}\biggl\lfloor\frac{b}{{\lambda}\log(1/{\lambda})}-
\frac{\mu}{{\lambda}\log^2(1/{\lambda})} \biggr\rfloor\biggr]\hspace*{-3.5pt}\biggr]
\subset\mathbb{Z},
\\
{\tilde Z}^{{\lambda},\mu}_t(J)
&=& 1 - \frac{\log(1+\#\{k\in
J_{{\lambda},\mu}, \eta^{\lambda}_{t\log
(1/{\lambda}
)}(k)=0\})}
{\log(1+\#(J_{{\lambda},\mu}))}.\nonumber
\end{eqnarray}
Observe that ${\tilde Z}^{{\lambda},\mu}_t(J)=1$ if and only if all
the sites of
$J_{{\lambda},\mu}$
are occupied at time $t\log(1/{\lambda})$. The quantity ${\tilde
Z}^{{\lambda},\mu}_t(J)$
is a function of the density of vacant clusters in the (rescaled) zone
$J$. Under some exchangeability properties, it should be closely related
to the size of occupied clusters in that zone, that is, to $Z^{\lambda}_t(x)$
for $x\in J$.

For $x\in(-A,A)$, ${\lambda}\in(0,1)$ and $\mu\in(0,1]$, we introduce
%
%
\begin{eqnarray} \label{xla}
x_{{\lambda},\mu} &=& \biggl[\hspace*{-3.5pt}\biggl[\biggl\lfloor\frac{x}{{\lambda}\log
(1/{\lambda})}-\frac{\mu
}{{\lambda}\log
^2(1/{\lambda})}
\biggr\rfloor+1,\nonumber\\
&&\hspace*{9.9pt} \biggl\lfloor\frac{x}{{\lambda}\log(1/{\lambda})}+
\frac{\mu}{{\lambda}\log^2(1/{\lambda})} \biggr\rfloor-1
\biggr]\hspace*{-3.5pt}\biggr]\subset\mathbb{Z},\\
{\tilde H}^{{\lambda},\mu}_t(x) &=& \frac{\log(1+\#\{k\in x_{{\lambda
},\mu}, \eta^{\lambda}_{t\log
(1/{\lambda}
)}(k)=0\})}
{\log(1+\#(x_{{\lambda},\mu}))}. \nonumber
\end{eqnarray}
Here, again, ${\tilde H}^{{\lambda},\mu}_t(x)=0$ if and only if all
the sites of
$x_{{\lambda},\mu}$
are occupied at time $t\log(1/{\lambda})$. Assume that a microscopic
fire starts
at some $x$. The process
${\tilde H}^{{\lambda},\mu}_t(x)$ will then allow us to quantify the
duration for which
this fire will be in effect.

Observe that we always have $\log(1+\#(x_{{\lambda},\mu}))\sim\log
(1+\#(J_{{\lambda},\mu}))
\sim\log(1/{\lambda})$ as ${\lambda}\to0$.
Also, observe that if ${\tilde Z}^{{\lambda},\mu}_t(J)=z$, then there are
$(1+\#(J_{{\lambda},\mu}))^{1-z}-1 \simeq{\lambda}^{z-1}$
vacant sites in $J_{{\lambda},\mu}$ at time $t\log(1/{\lambda})$.
In the same way, ${\tilde H}^{{\lambda},\mu}_t(x)=h$ says that there are
$(1+\#(x_{{\lambda},\mu}))^{h}-1 \simeq{\lambda}^{-h}$
vacant sites in $x_{{\lambda},\mu}$ at time $t\log(1/{\lambda})$.

We work conditionally on $M$.
We denote by ${\mathbb{P}}_M$ the conditional\break probability given $M$.
We recall that,
conditionally on $M$, $(Z_t(x),D_t(x),\break H_t(x))_{t\in[0,T],x\in[-A,A]}$
is deterministic.
We set $n =M([0,T] \times[-A;A])$, which is a.s. finite.
We set $T_0=0$ and consider the marks
$(X_q,T_q)_{1\leq q \leq n}$
of $M$, ordered in such a way that $T_0 < T_1< \cdots< T_n<T$.

We set ${\mathcal B}_0=\varnothing$ and for $q=1,\ldots,n$,
we consider ${\mathcal B}_q=\{X_1,\ldots,X_q\}$, as well as the set
${\mathcal C}_q$ of connected components of $(-A,A)\setminus{\mathcal
B}_q$ (sometimes
referred to as \textit{cells}).

Observe that, by construction, we have, for $c\in{\mathcal C}_q$ and
$x,y\in c$,
$Z_t(x)=Z_t(y)$ for all $t\in[0,T_{q+1})$. Thus, we can introduce
$Z_t(c)$.

We consider ${\lambda}_\mu>0$ (which depends on $M$)
such that for all ${\lambda}\in(0,{\lambda}_\mu)$, we have
$(X_i)_{{\lambda},\mu}\ne\varnothing$ and $(X_i)_{{\lambda},\mu
}\cap(X_j)_{{\lambda},\mu}=\varnothing$
for all $i\ne j$ with $i,j \in\{1,\ldots,n\}$.

We then observe that for ${\lambda}\in(0,{\lambda}_\mu)$ and for
each $q=0,\ldots,n$,
$\{x_{{\lambda},\mu},x\in{\mathcal B}_q\}\cup\{c_{{\lambda},\mu},
c \in{\mathcal C}_q \}$ is a partition
of $\llbracket-{\tilde A}_{{\lambda},\mu},{\tilde A}_{{\lambda},\mu
}\rrbracket$, where ${\tilde A}_{{\lambda},\mu}
=\lfloor
A/({\lambda}\log(1/{\lambda})) -
\mu/({\lambda}\log^2(1/{\lambda}))\rfloor$.

With our coupling, for the $({\lambda},A)$-FFP $(\eta^{\lambda
}_t)_{t\geq0}$,
for each $i= 1,\ldots,n$, a fire starts at the site
$\lfloor X_i/({\lambda}\log(1/{\lambda})) \rfloor$ at time
$T_i\log(1/{\lambda})$ and this describes all of the fires during
$[0,T\log
(1/{\lambda})]$.

The lemma below shows some exchangeability properties inside cells
[connected components of $(-A,A)\setminus{\mathcal B}_q$]. This will
allow us to prove that for $c$ a cell and $x\in c$, the size of
the occupied cluster around $x$
[described by $Z^{\lambda}(x)$] is closely related to the global
density of occupied
clusters in $c$ [described by ${\tilde Z}^{{\lambda},\mu}(c)$].
\begin{lem}\label{exch}
For ${\lambda}\in(0,1)$ and $\mu\in(0,1]$, set ${\mathcal
E}^{{\lambda},\mu}_0=\Omega$,
and for
$q=1,\ldots,n$, consider the event [recalling Definition \ref{dflaffpA}
and (\ref{xxbb})]
\begin{eqnarray*}
{\mathcal E}^{{\lambda},\mu}_q &=& \bigl\{\forall i=1,\ldots,q, \forall c \in
{\mathcal C}_{i},
\mbox{ either } c_{{\lambda},\mu}\subset C^{\lambda}_{T_i\log
(1/{\lambda})-}(X_i) \\
&&\hspace*{4.5pt}\mbox{or } \eta^{\lambda}_{T_i\log(1/{\lambda})-}(k)=0 \mbox{
for some }
\max c_{{\lambda},\mu}< k < \min C^{\lambda}_{T_i\log(1/{\lambda
})-}(X_i) \\
&&\hspace*{4.5pt}\mbox{or } \eta^{\lambda}_{T_i\log(1/{\lambda})-}(k)=0 \mbox{
for some }
\max C^{\lambda}_{T_i\log(1/{\lambda})-}(X_i)<k<\min c_{{\lambda
},\mu}\bigr\}.
\end{eqnarray*}
Conditionally on $M$ and ${\mathcal E}^{{\lambda},\mu}_q$, for all
$c\in{\mathcal C}_q$, the
random variables\break $(\eta^{\lambda}_{T_q\log(1/{\lambda})}(k))_{k\in
c_{{\lambda},\mu}}$ are
exchangeable.
\end{lem}
\begin{pf}
Let $c \in{\mathcal C}_q$, let $\sigma$ be a permutation of
$c_{{\lambda},\mu}$
and set, for simplicity, $\sigma(i)=i$ for $i\in I^{\lambda
}_A\setminus
c_{{\lambda},\mu}
$ [recall (\ref{alaila})].

Consider the $({\lambda},A)$-FFP process
$(\eta_t^{\lambda})_{t\geq0}$ constructed with $M$
and the family of Poisson processes $(N(i))_{i\in I^{\lambda}_A}$.
Also, consider the $({\lambda},A)$-FFP process
$(\tilde\eta_t^{\lambda})_{t\geq0}$ constructed with $M$
and the family of Poisson processes $({\tilde N}(i))_{i\in I^{\lambda}_A}$
defined by ${\tilde N}(i)=N(\sigma(i))$.

Observe that ${\mathcal E}^{{\lambda},\mu}_{k+1}\subset{\mathcal
E}^{{\lambda},\mu}_k$.
For all $k=0,\ldots,q$, $c\subset c_k$ for some $c_k\in{\mathcal C}_k$.
We will prove the following claims by induction on $k=0,\ldots,q$:

{\smallskipamount=0pt
\begin{longlist}
\item
if $\tilde{\mathcal E}^{{\lambda},\mu}_k$ is the same event as
${\mathcal E}^{{\lambda},\mu}_k$
corresponding to $(\tilde\eta^{\lambda}_t)_{t\geq0}$,
then $\tilde{\mathcal E}^{{\lambda},\mu}_k={\mathcal E}^{{\lambda
},\mu}_k$;

\item
on ${\mathcal E}^{{\lambda},\mu}_k$, for all $t\in[0,T_k\log
(1/{\lambda})]$,
$\tilde\eta_t^{\lambda}(i)= \eta_t^{\lambda}(\sigma(i))$ for all
$i\in I^{\lambda}_A$
[in particular, $\tilde\eta_t^{\lambda}(i)= \eta_t^{\lambda}(i)$
for all $i\notin c_{{\lambda},\mu}]$.
\end{longlist}}

Of course, (i) and (ii) with $k=q$ imply the lemma.
Indeed, let\break $\varphi\dvtx\{0,1\}^{\#(c_{{\lambda},\mu})}\mapsto{\mathbb{R}}$.
We have
\[
{\mathbb{E}}_M\bigl[{\mathbf{1}}_{{\mathcal E}^{{{\lambda},\mu
}}_q}\varphi\bigl(\bigl(\eta^{\lambda}_{T_q\log(1/{\lambda}
)}(i)\bigr)_{i\in
c_{{\lambda},\mu}} \bigr)\bigr]
={\mathbb{E}}_M\bigl[{\mathbf{1}}_{\tilde{\mathcal E}^{{{\lambda},\mu}}_q}
\varphi\bigl(\bigl(\tilde\eta^{\lambda}_{T_q\log(1/{\lambda})}(i)\bigr)_{i\in
c_{{\lambda},\mu}} \bigr)\bigr].
\]
Using (i) and (ii), we then deduce that
\[
{\mathbb{E}}_M\bigl[{\mathbf{1}}_{{\mathcal E}^{{{\lambda},\mu
}}_q}\varphi\bigl(\bigl(\eta^{\lambda}_{T_q\log(1/{\lambda}
)}(i)\bigr)_{i\in
c_{{\lambda},\mu}} \bigr)\bigr]
={\mathbb{E}}_M\bigl[{\mathbf{1}}_{{\mathcal E}^{{{\lambda},\mu}}_q}
\varphi\bigl(\bigl(\eta^{\lambda}_{T_q\log(1/{\lambda})}(\sigma(i))\bigr)_{i\in
c_{{\lambda},\mu}} \bigr)\bigr],
\]
which proves the lemma.

First, (i) and (ii) with $k=0$ are obviously satisfied.
Assume, now, that for some $k\in\{0,\ldots,q-1\}$, we have (i) and (ii).
Then, on ${\mathcal E}^{{\lambda},\mu}_k$,
for all $t\in[0,T_{k+1}\log(1/{\lambda}))$,
$\tilde\eta_t^{\lambda}(i)= \eta_t^{\lambda}(\sigma(i))$ for all
$i\in I^{\lambda}_A$.
Indeed, they are equal on $[0,T_{k}\log(1/{\lambda})]$, by assumption,
and they use the same
Poisson process ${\tilde N}(i)=N(\sigma(i))$ on the time interval
$[T_{k}\log(1/{\lambda}),T_{k+1}\log(1/{\lambda}))$).

We now check that ${\mathcal E}_{k+1}^{{\lambda},\mu}=\tilde
{\mathcal E}_{k+1}^{{\lambda},\mu}$. We know that
${\mathcal E}_{k}^{{\lambda},\mu}=\tilde{\mathcal E}_{k}^{{\lambda
},\mu}$ and the additional condition
[at time $T_{k+1}\log(1/{\lambda})-$] concerns:

$\bullet$ sites outside $c_{{\lambda},\mu}$, for which the values of
$\eta^{\lambda}
$ and
$\tilde\eta^{\lambda}$ at time $T_{k+1}\log(1/{\lambda})-$ are the same;

$\bullet$ the event $c_{{\lambda},\mu}\subset C^{\lambda
}_{T_{k+1}\log(1/{\lambda})-}$,
which is the same for $\eta^{\lambda}$ and $\tilde\eta^{\lambda}$
(it can be realized
only if there are no vacant sites in $c_{{\lambda},\mu}$, which
occurs, or not,
simultaneously for $\eta^{\lambda}$ and $\tilde\eta^{\lambda}$).

We now conclude that (ii) remains true at time $T_{k+1}\log(1/{\lambda
})$ since
the zone subject to fire either:

$\bullet$ is disjoint with $c_{{\lambda},\mu}$ so that the values of
$\eta^{\lambda},
\tilde\eta^{\lambda}$ are left invariant in $c_{{\lambda},\mu}$,
while they are modified
in the same way outside $c_{{\lambda},\mu}$;
or

$\bullet$ contains the whole zone $c_{{\lambda},\mu}$, which is thus
destroyed simultaneously for $\eta^{\lambda}$ and $\tilde\eta
^{\lambda}$, and the
values of
$\eta^{\lambda},\tilde\eta^{\lambda}$ are modified
in the same way outside $c_{{\lambda},\mu}$.
\end{pf}

The next lemma shows, in some sense, that
if a cell is \textit{almost} completely occupied at time $t$, then
it will be \textit{really} completely occupied at time $t+$;
and, if the effect of a microscopic fire is \textit{almost} ended at time
$t$, then
it will be \textit{really} ended at time $t+$.
\begin{lem}\label{undoncplein}
Let $\mu\in(0,1]$. Consider $k\in\{0,\ldots,n\}$, $c\in{\mathcal
C}_k$, $x\in
{\mathcal B}
_k$ and $t\in[T_k,T_{k+1})$.

{\smallskipamount=0pt
\begin{longlist}
\item Assume that for all ${\varepsilon}>0$, $\lim_{{\lambda}\to
0}{\mathbb{P}}
_M({\tilde Z}^{{\lambda},\mu}_t(c)<1-{\varepsilon})=0$.
Then, for all $s\in(t,T_{k+1})$, $\lim_{{\lambda}\to0}{\mathbb
{P}}_M({\tilde Z}^{{\lambda},\mu}
_{s}(c)=1)=1$.

\item Assume that for all ${\varepsilon}>0$, $\lim_{{\lambda}\to
0}{\mathbb{P}}
_M({\tilde H}^{{\lambda},\mu}_t(x)>{\varepsilon})=0$.
Then, for all $s\in(t,T_{k+1})$, $\lim_{{\lambda}\to0}{\mathbb
{P}}_M({\tilde H}^{{\lambda},\mu}
_{s}(x)=0)=1$.
\end{longlist}}
\end{lem}
\begin{pf} The proofs of (i) and (ii) are similar. Let us, for example,
prove (i). Thus, let $T_k\leq t< t+{\varepsilon}=s<T_{k+1}$. We start with
\[
{\mathbb{P}}_M\bigl({\tilde Z}^{{\lambda},\mu}_{t+{\varepsilon}}(c)=1 \bigr)
\geq
{\mathbb{P}}_M\bigl({\tilde Z}^{{\lambda},\mu}_{t+{\varepsilon}}(c)=1
{\mid}{\tilde Z}^{{\lambda},\mu}
_{t}(c)>1-{\varepsilon}/2\bigr) {\mathbb{P}}_M\bigl( {\tilde Z}
^{{\lambda},\mu}_{t}(c)>1-{\varepsilon}/2\bigr),
\]
so that it suffices to check that
$\lim_{{\lambda}\to0} {\mathbb{P}}_M({\tilde Z}^{{\lambda},\mu
}_{t+{\varepsilon}}(c)=1 {\mid}
{\tilde Z}^{{\lambda},\mu}
_{t}(c)>1-{\varepsilon}/2)=1$.
Let $v_t^{{\lambda},\mu}$
denote the number of vacant sites in $c_{{\lambda},\mu}$
(for $\eta^{\lambda}_{t\log(1/{\lambda})}$).
Then ${\tilde Z}^{{\lambda},\mu}_{t+{\varepsilon}}(c)=1$ is
equivalent to $v^{{\lambda},\mu}
_{t+{\varepsilon}}=0$ and
one can easily check that
${\tilde Z}^{{\lambda},\mu}_t(c)>1-{\varepsilon}/2$ implies that
$v^{{\lambda},\mu}_t \leq(1+\#
(c_{{\lambda},\mu}))^{{\varepsilon}
/2}\leq
(1+2A/({\lambda}\log(1/{\lambda})))^{{\varepsilon}/2}$.

Since $M((t,s]\times[-A,A])=0$ by assumption,
we deduce that $M^{\lambda}_{s\log(1/{\lambda})}(i)=M^{\lambda
}_{t\log(1/{\lambda})}(i)$
for all\vspace*{1pt} $i\in I^{\lambda}_A$: no fire starts during $(t\log
(1/{\lambda}),s\log
(1/{\lambda})]$.
Hence, each occupied site at time $t\log(1/{\lambda})$
remains occupied at time $s\log(1/{\lambda})$ and each vacant site at time
$t\log(1/{\lambda})$
becomes occupied at time $s\log(1/{\lambda})$ with probability
$1-e^{(t-s)\log(1/{\lambda})}=1-{\lambda}^{\varepsilon}$.
Thus,
\[
{\mathbb{P}}_M\bigl({\tilde Z}^{{\lambda},\mu}_{t+{\varepsilon}}(c)=1
{\mid}{\tilde Z}^{{\lambda},\mu}
_{t}(c)>1-{\varepsilon}/2\bigr)
\geq(1-{\lambda}^{\varepsilon})^{(1+2A/({\lambda}\log(1/{\lambda
})))^{{\varepsilon}/2}},
\]
which tends to $1$ as ${\lambda}\to0$.
\end{pf}

We end our preliminaries with a last lemma, which deals with estimates
concerning the time needed to occupy vacant zones.
\begin{lem}\label{binomiale}
Let $\mu\in(0,1]$. Let $(\zeta_0^{\lambda}(i))_{i\in I^{\lambda
}_A}\in\{0,1\}
^{I^{\lambda}_A}$ and
consider a family of i.i.d. Poisson processes
$(P^{\lambda}_t(i))_{t\geq0,i\in I^{\lambda}_A}$, with rate $\log
(1/{\lambda})$, independent
of $\zeta^{\lambda}_0$.
Set $\zeta^{\lambda}_t(i)=\min(\zeta_0^{\lambda}(i)+P^{\lambda}_t(i),1)$.

\begin{enumerate}
\item Let $J=(a,b)\subset(-A,A)$ and $h\in[0,1]$.
Set $v_t^{{\lambda},\mu}=\#\{i\in J_{{\lambda},\mu}, \zeta
^{\lambda}_t(i)=0\}$.
Assume that
\[
\forall{\varepsilon}>0\qquad
{\mathbb{P}}\biggl( \biggl|\frac{\log(1+v^{{\lambda},\mu}_0)}{\log(1+\#
(J_{{\lambda},\mu}))} - h
\biggr|\geq{\varepsilon}\biggr)=0.
\]

\begin{enumerate}[(a)]
\item[(a)]
Then, for all $T>0$ and ${\varepsilon}>0$,
\[
\lim_{{\lambda}\to0}
{\mathbb{P}}\biggl(\sup_{[0,T]} \biggl|\frac{\log(1+v^{{\lambda},\mu
}_t)}{\log(1+\#(J_{{\lambda},\mu}))} - (h-t)_+
\biggr|\geq
{\varepsilon}\biggr)=0.
\]
\item[(b)] If the family $(\zeta^{\lambda}_0(i))_{i\in J_{{\lambda},\mu}}$
is exchangeable, then, for all $x\in J$,
$T>0$ and ${\varepsilon}>0$,
\[
\lim_{{\lambda}\to0}
{\mathbb{P}}\biggl(\sup_{[0,T]} \biggl|\frac{\log(1+\# (G^{\lambda
}_t(x)))}{\log(1/{\lambda})}
- \bigl(1-(h-t)_+\bigr) \biggr|
\geq{\varepsilon}\biggr)=0,
\]
where $G^{\lambda}_t(x)$ is the connected component of occupied
sites around $\lfloor x/{\lambda}\log(1/{\lambda})\rfloor$
in $\zeta^{\lambda}_t$.
\end{enumerate}

\item Let $x\in(-A,A)$ and $h\in[0,1]$.
Set $v_t^{{\lambda},\mu}=\#\{i\in x_{{\lambda},\mu}, \zeta
^{\lambda}_t(i)=0\}$.
Assume that
\[
\forall{\varepsilon}>0\qquad
{\mathbb{P}}\biggl( \biggl|\frac{\log(1+v^{{\lambda},\mu}_0)}{\log(1+\#
(x_{{\lambda},\mu}))} - h \biggr|
\geq{\varepsilon}\biggr)=0.
\]
Then, for all $T>0$ and ${\varepsilon}>0$,
\[
\lim_{{\lambda}\to0}
{\mathbb{P}}\biggl(\sup_{[0,T]} \biggl|\frac{\log(1+v^{{\lambda},\mu
}_t)}{\log(1+\#(x_{{\lambda},\mu}))}
- (h-t)_+ \biggr|\geq
{\varepsilon}\biggr)=0.
\]
\end{enumerate}
%
\end{lem}
\begin{pf}
The proof of part 2 is the same as that of 1(a)
because $\log(1+\#(J_{{\lambda},\mu}))\sim\log(1+\#(x_{{\lambda
},\mu}))\sim\log(1/{\lambda})$
as ${\lambda}\to0$.
Thus, we only prove 1 and everywhere replace
$\log(1+\#(x_{{\lambda},\mu}))$ by $\log(1/{\lambda})$ without difficulty.
By assumption, for all ${\varepsilon}>0$, we have $\lim_{{\lambda
}\to0}
{\mathbb{P}}(v^{{\lambda},\mu}_0 \in({\lambda}^{{\varepsilon
}-h}-1,{\lambda}^{-{\varepsilon
}-h})) =1$.
We define $h_t=(h-t)_+$, $V^{{\lambda},\mu}_t=\log(1+v^{{\lambda
},\mu}_t)/\log(1/{\lambda})$
and, finally, $\Gamma^{\lambda}_t=\log(1+\#(G^{\lambda}_t(x)))/\log
(1/{\lambda})$.\vspace*{8pt}

\textit{Step} 1. Let $t\geq0$ be fixed.
We first show that for all ${\varepsilon}>0$,
$\lim_{{\lambda}\to0}{\mathbb{P}}(|V^{{\lambda},\mu}_t - h_t
|\geq{\varepsilon})=0$.
Conditionally on $v^{{\lambda},\mu}_0$, the random variable
$v^{{\lambda},\mu}_t$ follows a
binomial distribution $B(v_0^{{\lambda},\mu},{\lambda}^t)$ because
each vacant site at time
$0$ remains vacant with probability $e^{-t\log(1/{\lambda})}={\lambda}^t$.\vspace*{8pt}

\textit{Case $h_t>0$.} Let ${\varepsilon}\in(0,h_t)$.
We have to prove that ${\mathbb{P}}(v^{{\lambda},\mu}_t\in({\lambda
}^{{\varepsilon}-h_t}$,\break ${\lambda}
^{-{\varepsilon}-h_t})
)\to1$. We know that $\lim_{{\lambda}\to0} {\mathbb{P}}(
v^{{\lambda},\mu}_0 \in({\lambda}^{{\varepsilon}/2-h},{\lambda
}^{-{\varepsilon}/2-h}) )=1$.
The Bienaym\'e--Chebyshev inequality
implies that
\begin{eqnarray*}
&&P[|v^{{\lambda},\mu}_t-v_0^{\lambda}{\lambda}^t|\leq
(v_0^{{\lambda},\mu}{\lambda}^t)^{2/3}{\mid}
v_0^{{\lambda},\mu}
\in({\lambda}^{{\varepsilon}/2-h},{\lambda}^{-{\varepsilon
}/2-h})]\\
&&\qquad\geq1- {\mathbb{E}}[v_0^{{\lambda},\mu}{\lambda}^t (1-{\lambda
}^t) (v_0^{{\lambda},\mu}{\lambda}
^t)^{-4/3}{\mid}
v_0^{{\lambda},\mu}
\in({\lambda}^{{\varepsilon}/2-h},{\lambda}^{-{\varepsilon}/2-h})]
\\
&&\qquad\geq1 - {\mathbb{E}}[(v_0^{{\lambda},\mu}{\lambda}^t)^{-1/3}{\mid
}v_0^{{\lambda},\mu}
\in({\lambda}^{{\varepsilon}/2-h},{\lambda}^{-{\varepsilon}/2-h})]\\
&&\qquad\geq1- ({\lambda}^{{\varepsilon}/2-h+t})^{-1/3},
\end{eqnarray*}
which tends to $1$ since $h_t=h-t>{\varepsilon}$.

However, the events
\[
|v^{{\lambda},\mu}_t-v_0^{{\lambda},\mu
}{\lambda}^t|\leq(v_0^{{\lambda},\mu}{\lambda}
^t)^{2/3}
\quad\mbox{and}\quad
v_0^{{\lambda},\mu}\in({\lambda}^{{\varepsilon}/2-h}, {\lambda
}^{-{\varepsilon}/2-h})
\]
imply that
$v^{{\lambda},\mu}_t \in
({\lambda}^{{\varepsilon}/2-h_t}-({\lambda}^{-{\varepsilon
}/2-h_t})^{2/3},{\lambda}
^{-{\varepsilon}/2-h_t}+({\lambda}^{-{\varepsilon}
/2-h_t})^{2/3})
\subset({\lambda}^{{\varepsilon}-h_t}$,\break ${\lambda}^{-{\varepsilon
}-h_t})$ for
${\lambda}$ small enough, hence the result.\vspace*{8pt}

\textit{Case $h_t=0$.} We have to show that for all ${\varepsilon}>0$,
$\lim_{{\lambda}\to0}{\mathbb{P}}(v^{{\lambda},\mu}_t>{\lambda
}^{-{\varepsilon}})=0$, and it
suffices to check that
$\lim_{{\lambda}\to0}{\mathbb{P}}(v^{{\lambda},\mu}_t>{\lambda
}^{-{\varepsilon}} {\mid}v^{{\lambda},\mu}
_0<{\lambda}^{-{\varepsilon}/2-h})=0$.
However,
\begin{eqnarray*}
&&{\mathbb{P}}(v^{{\lambda},\mu}_t>{\lambda}^{-{\varepsilon}} {\mid
}v^{{\lambda},\mu}_0<{\lambda}
^{-{\varepsilon}/2-h})\\
&&\qquad\leq{\lambda}^{\varepsilon}
{\mathbb{E}}[v^{{\lambda},\mu}_t {\mid}v^{{\lambda},\mu
}_0<{\lambda}^{-{\varepsilon}/2-h}] = {\lambda}
^{\varepsilon}
{\mathbb{E}}[ v^{{\lambda},\mu}_0 {\lambda}^t {\mid}v^{{\lambda
},\mu}_0<{\lambda}^{-{\varepsilon}/2-h}]\\
&&\qquad\leq{\lambda}^{{\varepsilon}+t}{\lambda}^{-{\varepsilon}/2-h}
={\lambda}
^{{\varepsilon}/2+t-h},
\end{eqnarray*}
which tends to $0$ since, by assumption, $t-h\geq0$.\vspace*{8pt}

\textit{Step} 2. We now prove that for all ${\varepsilon}>0$, $\lim
_{{\lambda}
\to0} {\mathbb{P}}
(|\Gamma^{\lambda}_t-(1-h_t)|\geq{\varepsilon})=0$.
It suffices to check that $\lim_{{\lambda}\to0} {\mathbb{P}}
( \#(G^{\lambda}_t(x)) \in({\lambda}^{{\varepsilon
}+h_t-1}-1,{\lambda}
^{-{\varepsilon}+h_t-1}))=1$.
However, we know from Step 1 that there are approximately
$(1/{\lambda})^{h_t}$ vacant sites in $J_{{\lambda},\mu}$, and
$\#(J_{{\lambda},\mu})\simeq(1/{\lambda}\log(1/{\lambda}))$.
We also know that
the family $(\zeta^{\lambda}_t(i))_{i\in J_{{\lambda},\mu}}$ is
exchangeable so that
the vacant sites are
uniformly distributed in $J_{{\lambda},\mu}$ (this statement is slightly
misleading: there
cannot be two vacant sites at the same place).
We conclude that $\#(G^{\lambda}_t(x)) \simeq(1/{\lambda}\log
(1/{\lambda}))/\break
(1/{\lambda})^{h_t} \simeq{\lambda}^{h_t-1}$.
This can be done rigorously without difficulty.\vspace*{8pt}

\textit{Step} 3. We now prove 1(a), which relies on Step 1 and an ad hoc
version of Dini's theorem. Let ${\varepsilon}>0$. Consider
a subdivision $0=t_0<t_1<\cdots<t_l=T$ with $t_{i+1}-t_i<{\varepsilon}/2$.
Using Step 1, we have
$\lim_{{\lambda}\to0}{\mathbb{P}}[\max_{i=0,\ldots,l}|V^{{\lambda
},\mu}_{t_i} -
(h-t_i)_+|>{\varepsilon}/2]=0$.

Now, observe that
$t\mapsto V_t^{{\lambda},\mu}$ and $t \mapsto(h-t)_+$ are a.s. nonincreasing
and that $t \mapsto(h-t)_+$ is Lipschitz continuous with Lipschitz
constant $1$.

We deduce that ${\sup_{[0,T]}} |V^{{\lambda},\mu}_t-(h-t)_+|
\leq{\varepsilon}/2+ \max_{i=0,\ldots,l}\{|V^{{\lambda},\mu}_{t_i} -
(h-t_i)_+|\}$.
Thus, ${\mathbb{P}}({\sup_{[0,T]}} |V^{{\lambda},\mu}_t-(h-t)_+|
>{\varepsilon})
\leq{\mathbb{P}}[{\max_{i=0,\ldots,l}}|V^{{\lambda},\mu}_{t_i} -
(h-t_i)_+|>{\varepsilon}/2]$,
which completes the proof of 1(a).\vspace*{8pt}

\textit{Step} 4. Point 1(b) is deduced from Step 2 exactly as
point 1(a) was deduced from Step 1, using the fact that $t\mapsto
\Gamma
^{\lambda}_t$
and $t \mapsto1-h_t$ are a.s. nondecreasing.
\end{pf}

We may now finally tackle the following proof.
\begin{pf*}{Proof of Proposition \ref{heart}}
For $x\in(-A,A)$ and $t\geq0$, we introduce
$Z_t(x-)=\lim_{y\to x, y<x}Z_t(y)$
and $Z_t(x+)=\lim_{y\to x, y>x}Z_t(y)$, which represent the values
of $Z_t$ in the cells on the left and right of $x$.
If $x \in{\mathcal B}_n$, it is at the boundary of two cells
$c_-,c_+ \in{\mathcal C}_n$, and then $Z_t(x-)=Z_t(c_-)$ and
$Z_t(x+)=Z_t(c_+)$.

For $x\in{\mathcal B}_n$ and $t\geq0$, we set ${\tilde H}_t(x)=\max
(H_t(x),1-Z_t(x),1-Z_t(x-),1-Z_t(x+))$.
Observe that for the LFFP, $x$ is \textit{microscopic} (or
\textit{acts like a barrier}) if and only if ${\tilde H}_t(x)>0$ and, if
so, it will
remain microscopic during exactly $[t,t+{\tilde H}_t(x))$.
Note that, in fact, $Z_t(x)$ always equals either $Z_t(x-)$ or $Z_t(x+)$.

We consider the set of times
${\mathcal K}:=\{t\in\{0,T\}$: there exists
$x\in(-A,A)$ such that ${\tilde H}_t(x)=0$ but ${\tilde
H}_{t-{\varepsilon}}(x)>0$
for all ${\varepsilon}>0$ small enough$\}$.
By construction, we see that ${\mathcal K}\subset\{
1,T_i+1,T_i+Z_{T_i-}(X_i),i=1,\ldots,n\}
\subset\{1, T_i+1,T_i+(T_i-T_j), 0\leq j < i \leq n \} $.

We work conditionally on $M$, by induction on $q=0,\ldots,n$.
Consider the following assumption.

{\smallskipamount=0pt
\begin{longlist}
\item[$({\mathcal H}_q)$: (i)]
For all $0< \mu\leq1$, $c\in{\mathcal C}_q$ and
${\varepsilon}>0$,
$\lim_{{\lambda}\to0}{\mathbb{P}}_M(|{\tilde Z}^{{{\lambda},\mu}
}_{T_q}(c)-Z_{T_q}(c)|>{\varepsilon})=0$.


\item[(ii)] For all $x\in{\mathcal B}_q$, $0 < \mu\leq1$ and ${\varepsilon}>0$,
$\lim_{{\lambda}\to
0}{\mathbb{P}}_M(|{\tilde H}^{{\lambda},\mu}_{T_q}(x)-{\tilde
H}_{T_q}(x)|>{\varepsilon})=0$.
%

\item[(iii)] For all $0 < \mu\leq1$, $\lim_{{\lambda}\to0}{\mathbb
{P}}_M({\mathcal E}^{{\lambda},\mu}_q)=1$
(recall Lemma \ref{exch}).
\end{longlist}}

First, $({\mathcal H}_0)$ is obviously satisfied because $T_0=0$,
${\mathcal C}_0=(-A,A)$,
${\tilde Z}^{{\lambda},\mu}_{0}((-A$,\break $A))=0=Z_{0}((-A,A))$, ${\mathcal
B}_0=\varnothing$
and ${\mathcal E}^{{\lambda},\mu}_0=\Omega$.

The proposition will essentially be proven if we check that for
$q=0,\ldots,n-1$, $({\mathcal H}_{q})$ implies:

(a) for $c\in{\mathcal C}_q$, $0<\mu\leq1$ and ${\varepsilon}>0$,
$\lim_{{\lambda}\to0}{\mathbb{P}}_M(\sup_{[T_q,T_{q+1})}|{\tilde
Z}^{{\lambda},\mu}
_{t}(c)-Z_{t}(c)|>{\varepsilon})=0$;

(b) for $x \in(-A,A)\setminus{\mathcal B}_q$, ${\varepsilon}>0$,
$\lim_{{\lambda}\to0}{\mathbb{P}}_M({\sup
_{[T_q,T_{q+1})}}|Z^{\lambda}
_{t}(x)-Z_{t}(x)|>{\varepsilon})=0$;
%

(c) for $x\in{\mathcal B}_q$, $t\in[T_q,T_{q+1})$, $0< \mu\leq1$ and
${\varepsilon}
>0$, $\lim_{{\lambda}\to0}{\mathbb{P}}_M(|{\tilde H}^{{\lambda
},\mu}_{t}(x)-{\tilde H}
_{t}(x)|>{\varepsilon})$;

(d) for $x \in(-A,A)\setminus{\mathcal B}_q$,
$t\in(T_q,T_{q+1})\setminus{\mathcal K}$ and ${\varepsilon}>0$,
$\lim_{{\lambda}\to0}{\mathbb{P}}_M(\delta(D^{\lambda}_{t}(x)$,\break $D_{t}(x)
)>{\varepsilon})=0$;

(e) for $x \in(-A,A) \setminus{\mathcal B}_q$, ${\varepsilon}>0$,
$\lim_{{\lambda}\to0}{\mathbb{P}}_M(\int_{T_q}^{T_{q+1}}\delta
(D^{\lambda}_{t}(x) ,
D_{t}(x) ) \,dt
>{\varepsilon})=0$;

(f) $({\mathcal H}_{q+1})$ holds.

We thus assume $({\mathcal H}_{q})$ for some fixed $q\in\{0,\ldots,n-1\}$
and prove points (a)--(f). Below, we repeatedly use the fact that
on the time interval $[T_q,T_{q+1})$, there are no fires at all in $(-A,A)$
for the LFFP and no fires at all during
$[T_{q}\log(1/{\lambda}),T_{q+1}\log(1/{\lambda}))$ for the
${\lambda}$-FFP.

Set $\zeta^{\lambda}_0(i)=\eta^{\lambda}_{T_q\log(1/{\lambda})}(i)$
and consider the i.i.d. Poisson processes $P^{\lambda
}_t(i)=N_{(T_q+t)\log
(1/{\lambda})}(i)
-N_{T_q\log(1/{\lambda})}(i)$ with rate $\log(1/{\lambda})$. Then,
for $t\in
[T_q,T_{q+1})$,
$\eta^{\lambda}_{t\log(1/{\lambda})}(i)=\min(\zeta
_0(i)+P^{\lambda}_{t-T_q}(i),1)$.\vspace*{8pt}

\textit{Point} (a).
Let $0< \mu\leq1$. Let $c \in{\mathcal C}_q$. Observe that
$({\mathcal H}_{q})$(i)
says precisely that with
$h=1-Z_{T_q}(c)\in[0,1]$, $\log(1+\#\{k\in c_{{\lambda},\mu}, \zeta
^{\lambda}
_0(k)=0\})
/\log(1+\#(c_{{\lambda},\mu}))$ tends to $h$ in probability (for
${\mathbb{P}}_M$).
Applying part 1(a) of Lemma \ref{binomiale} (with $J=c$),
we get that $\sup_{[T_q,T_{q+1})} |1-{\tilde Z}^{{\lambda},\mu
}_t(c)-(h-(t-T_q))_+|$
tends to $0$ in probability (for ${\mathbb{P}}_M$). However,
for $t\in[T_q,T_{q+1})$,
we have $Z_t(c)=\min(Z_{T_q}(c)+(t-T_q),1)=\min(1-h+(t-T_q),1)
=1-(h-(t-T_q))_+$. Point (a) then follows.\vspace*{8pt}

\textit{Point} (b).
Now, let $x\in(-A,A)\setminus{\mathcal B}_q$. Then $x \in c$, for
some $c\in
{\mathcal C}_q$.
Due to Lemma \ref{exch}, we know that $(\zeta_0^{\lambda}(i))_{i \in
c_{{\lambda},\mu}}$
are exchangeable on ${\mathcal E}^{{\lambda},1}_q$.
The previous reasoning, using part 1(b)
of part 1(a) of Lemma \ref{binomiale}, shows that for all
${\varepsilon}>0$,
$\lim_{{\lambda}\to0}{\mathbb{P}}_M( {\mathcal E}^{{\lambda},1}_q
\cap\{
{\sup_{[T_q,T_{q+1})}}|Z^{\lambda}_{t}(x)-Z_{t}(x)|>{\varepsilon}\}
)=0$. Using
$({\mathcal H}_{q})$(iii) for $\mu=1$, we are done.\vspace*{8pt}

\textit{Point} (c). Let $0< \mu\leq1$. Let $x \in{\mathcal B}_q$ and
set $h={\tilde H}_{T_q}(x)$.
We know by $({\mathcal H}_q)$(ii) that ${\tilde H}_{T_q}^{{\lambda},\mu
}(x)$ tends to ${\tilde H}_{T_q}(x)=h$
in probability (for ${\mathbb{P}}_M$). Now, using part 2(a) of Lemma
\ref{binomiale},
we deduce that ${\sup_{[T_q,T_{q+1})}}|{\tilde H}^{{\lambda},\mu
}_t(x)-(h-(t-T_q))_+ |$
tends to $0$ in probability (for ${\mathbb{P}}_M$). We conclude by observing
that, by construction, ${\tilde H}_t(x)= (h-(t-T_q))_+$ for $t\in
[T_q,T_{q+1})$.\vspace*{8pt}

\textit{Point} (d). Let $x \in(-A,A)\setminus{\mathcal B}_q$
and $t \in(T_q,T_{q+1})\setminus{\mathcal K}$ be fixed.\vspace*{8pt}

\textit{Case $Z_t(x)<1$.} In this case, $D_t(x)=\{x\}$ so that
$\delta(D_t(x),D^{\lambda}_t(x))=|D^{\lambda}_t(x)|$. However,
from (\ref{dlambda}), (\ref{zlambda}), we get that $|D^{\lambda}_t(x)|
\leq
{\lambda}^{1-Z^{\lambda}_t(x)}
\log(1/{\lambda})$. Since we know from (b) that $Z^{\lambda}_t(x)$
goes to $Z_t(x)<1$
in probability (for ${\mathbb{P}}_M$), we easily deduce that
$|D^{\lambda}_t(x)|$
goes to $0$ in probability (for ${\mathbb{P}}_M$).\vspace*{8pt}

\textit{Case $Z_t(x)=1$.} In this case, $D_t(x)=[a,b]$ for some
$a,b\in{\mathcal B}_q\cup\{-A,A\}$. We assume that $-A<a<b<A$ for simplicity,
the other cases being treated in a similar way. We thus
have $Z_t(c)=1$ for all $c\in{\mathcal C}_q$
with $c\subset(a,b)$, ${\tilde H}_t(y)=0$ for all $y\in{\mathcal
B}_q\cap(a,b)$
and ${\tilde H}_t(a){\tilde H}_t(b)>0$.

On the one hand, we prove that for any ${\varepsilon}>0$,
$\lim_{{\lambda}\to0}{\mathbb{P}}_M(D^{\lambda}_t(x)\subset
[a-{\varepsilon
},b+{\varepsilon}])=1$.
Let us consider, for example, the left boundary $a$ and prove that
$\lim_{{\lambda}\to0}{\mathbb{P}}_M(D^{\lambda}_t(x)\subset
[a-{\varepsilon},A])=1$.

We have ${\tilde H}_t(a)=h_a>0$.
We deduce from (c) that
$\lim_{{\lambda}\to0} {\mathbb{P}}_M({\tilde H}^{{\lambda
},1}_t(a)\geq h_a/2)=1$,
which implies that there are vacant sites in $a_{{\lambda},1}$,
that is,\break
$\lim_{{\lambda}\to0} {\mathbb{P}}_M(\exists i\in a_{{\lambda},1},
\eta_{t\log(1/{\lambda}
)}(i)=0)=1$.
Recalling the definition of $a_{{\lambda},1}$ [see (\ref{xla})], we
see that
this implies that
$\lim_{{\lambda}\to0}{\mathbb{P}}_M(D^{\lambda}_t(x)\subset
[a-1/\log(1/{\lambda}),\break A])=1$,
hence
$\lim_{{\lambda}\to0}{\mathbb{P}}_M(D^{\lambda}_t(x)\subset
[a-{\varepsilon},A])=1$
for any
${\varepsilon}>0$.

%


On the other hand, we prove that $\lim_{{\lambda}\to0} {\mathbb
{P}}_M((a+1/\log
(1/{\lambda}
),b-1/\log(1/\break{\lambda}))
\subset D^{\lambda}_t(x))=1$.
Since $t\notin{\mathcal K}$, we deduce that there exists $s\in
(T_q,t)$ such that
$Z_s(c)=1$ for all $c\in{\mathcal C}_q$ with $c\subset(a,b)$
and ${\tilde H}_s(y)=0$ for all
$y\in{\mathcal B}_q\cap(a,b)$. We deduce from (a) that
for all $c\in{\mathcal C}_q$ with $c\subset(a,b)$,
$\lim_{{\lambda}\to0} {\mathbb{P}}_M({\tilde Z}^{{\lambda
},1}_s(c)>1-{\varepsilon})=0$,
whence, by Lemma
\ref{undoncplein}(i), $\lim_{{\lambda}\to0} {\mathbb{P}}_M({\tilde
Z}^{{\lambda},1}_t(c)=1)=1$.\break
Similarly, we deduce from (c) that
for all $y\in{\mathcal B}_q$ with $y\in(a,b)$,\break
$\lim_{{\lambda}\to0} {\mathbb{P}}_M({\tilde H}^{{\lambda
},1}_s(y)>{\varepsilon})=0$, whence,
by Lemma
\ref{undoncplein}(ii), $\lim_{{\lambda}\to0} {\mathbb
{P}}_M({\tilde H}^{{\lambda},1}_t(y)=0)=1$.
As a consequence, $\lim_{{\lambda}\to0} {\mathbb{P}}_M((a+1/\log
(1/{\lambda}),b-1/\log
(1/{\lambda}))
\subset D^{\lambda}_t(x))=1$.

This completes the proof of point (d).\vspace*{8pt}

\textit{Point} (e). Point (e) follows from (d). Indeed,
observe that $\delta(I,J) \leq2A$ for any intervals $I,J \subset(-A,A)$.
Thus, for $x\in(-A,A)\setminus{\mathcal B}_q$,
(d) implies
that for $t\in[T_q,T_{q+1})\setminus{\mathcal K}$,
$\lim_{{\lambda}\to0}{\mathbb{E}}_M(\delta(D^{\lambda}_t(x),D_t(x)))=0$.
Since ${\mathcal K}$ is now finite, we deduce
from Lebesgue's dominated convergence theorem that
$\lim_{{\lambda}\to0}\int_{T_q}^{T_{q+1}} {\mathbb{E}}_M(\delta
(D^{\lambda}_t(x),D_t(x)))\,dt=0$,
from which (e) follows.\vspace*{8pt}

\textit{Point} (f). Here, we show that $({\mathcal H}_{q+1})$ holds.
We set $z:=Z_{T_{q+1}-}(X_{q+1})$ and separately treat the cases $z\in(0,1)$
and $z=1$. We a.s. never have $z=0$ because
$Z_{T_{q+1}-}(X_{q+1})=\min(Z_{T_q}(X_{q+1})+(T_{q+1}-T_q),1)$
with $Z_{T_q}(X_{q+1})\geq0$
and $T_{q+1}>T_q$.\vspace*{8pt}

\textit{Case} $z\in(0,1)$. We fix $\mu\in(0,1]$. In that case,
$D_{T_{q+1}-}(X_{q+1})=\{X_{q+1}\}$ and for all
$c \in{\mathcal C}_{q+1}$ (thus $c\subset\tilde c$ for some $\tilde
c\in{\mathcal C}_q$),
$Z_{T_{q+1}}(c)=Z_{T_{q+1}-}(c)$.
We have ${\tilde H}_{T_{q+1}}(X_{q+1})=\max(z,1-z)$ and
for all $x\in{\mathcal B}_q$, ${\tilde H}_{T_{q+1}}(x)={\tilde
H}_{T_{q+1}-}(x)$.
Consider the event
$\Omega^{\lambda}_\alpha=\{Z^{\lambda}_{T_{q+1}-}(X_{q+1})\leq
z+\alpha\}$
for some $\alpha\in(0,1-z)$. Point (b) implies that
$\lim_{{\lambda}\to0}{\mathbb{P}}_M(\Omega^{\lambda}_\alpha)=1$ (because
$X_{q+1}\notin{\mathcal B}_q$).

\begin{itemize}
\item
On $\Omega^{\lambda}_\alpha$,
we have $\#(C^{\lambda}_{T_{q+1}\log(1/{\lambda})-}(X_{q+1}))\leq
(1/{\lambda}
)^{z+\alpha}$
[see (\ref{zlambda})].
Since $z+\alpha<1$, we deduce
that on $\Omega^{\lambda}_\alpha$, we have
$\#(C^{\lambda}_{T_{q+1}\log(1/{\lambda})-}(X_{q+1}))< \mu
/\break (2{\lambda}\log^2(1/{\lambda}))$
(for all $\mu$, provided that ${\lambda}>0$ is small enough).
Thus,
on $\Omega^{\lambda}_\alpha$, for all $c\in{\mathcal C}_{q+1}$,
there is a vacant site (strictly)
between $c_{{\lambda},\mu}$ and $C^{\lambda}_{T_{q+1}\log
(1/{\lambda})-}(X_{q+1})$. Hence,
${\mathcal E}^{{\lambda},\mu}_q\cap\Omega^{\lambda}_\alpha\subset
{\mathcal E}^{{\lambda},\mu}_{q+1}$. Using
$({\mathcal H}_q)$(iii),
we deduce that $\lim_{{\lambda}\to0} {\mathbb{P}}_M({\mathcal
E}^{{\lambda},\mu}_{q+1})=1$.

\item This also implies that on $\Omega^{\lambda}_\alpha$,
for all $c\in{\mathcal C}_{q+1}$, we have
${\tilde Z}^{{\lambda},\mu}_{T_{q+1}}(c)={\tilde Z}^{{\lambda},\mu
}_{T_{q+1}-}(c)$ and thus point (a)
and $\lim_{{\lambda}\to0}{\mathbb{P}}_M(\Omega^{\lambda}_\alpha)=1$
imply that $\lim_{{\lambda}\to0}{\mathbb{P}}_M(|{\tilde
Z}^{{\lambda},\mu}
_{T_{q+1}}(c)-Z_{T_{q+1}}(c)|\geq{\varepsilon})=0$
for all ${\varepsilon}>0$.

\item For $x\in{\mathcal B}_{q+1}\setminus\{X_{q+1}\}={\mathcal
B}_q$, still on
$\Omega
^{\lambda}_\alpha$,
we also have ${\tilde H}^{{\lambda},\mu}_{T_{q+1}}(x)={\tilde
H}^{{\lambda},\mu}_{T_{q+1}-}(x)$, thus
point (c)
allows us to conclude that $({\mathcal H}_{q+1})$(ii) holds for those points
$x$.

We now show that
$\lim_{{\lambda}\to0}{\mathbb{P}}_M(|{\tilde H}^{{\lambda},\mu
}_{T_{q+1}}(X_{q+1})-{\tilde H}
_{T_{q+1}}(X_{q+1})|\geq{\varepsilon})=0$
for all ${\varepsilon}>0$, which implies
that $({\mathcal H}_{q+1})$(ii) holds for $x=X_{q+1}$.
Recall that ${\tilde H}_{T_{q+1}}(X_{q+1})=\max(z,1-z)$.
Consider $c\in{\mathcal C}_q$ such that $X_{q+1}\in c$
and denote by $v^{{\lambda},\mu}_t$ the number of vacant sites in
$x_{{\lambda},\mu}$ at time
$t\log(1/{\lambda})$.
Point (a) implies that at time $T_{q+1}\log(1/{\lambda})-$,
there are around $(1/{\lambda})^{1-z}$ vacant sites in $c_{{\lambda
},\mu}$. Thus,
by exchangeability of the family $(\eta^{\lambda}_{T_{q+1}\log
(1/{\lambda}
)-}(i))_{i\in c_{{\lambda},\mu}}$
(on the event ${\mathcal E}^{{{\lambda},\mu}}_{q}$, see Lemma \ref{exch}),
since $x_{{\lambda},\mu}\subset c_{{\lambda},\mu}$ and $\#
(x_{{\lambda},\mu})/\#(c_{{\lambda},\mu})\simeq
1/\log
(1/{\lambda})$, 
we deduce that $v^{{\lambda},\mu}_{T_{q+1}-} \simeq(1/{\lambda
})^{1-z}/\log(1/{\lambda})
\simeq(1/{\lambda})^{1-z}$ on ${\mathcal E}^{{\lambda},\mu}_q$.
On the other hand, recalling (\ref{zlambda}), we have
$\#(C^{\lambda}_{T_{q+1}\log(1/{\lambda})-}(X_{q+1})) \simeq
(1/{\lambda})^z$. At time
$T_{q+1}\log(1/{\lambda})$, this component is destroyed. Thus,
still on ${\mathcal E}^{{\lambda},\mu}_q$,
$v^{{\lambda},\mu}_{T_{q+1}} = v^{{\lambda},\mu}_{T_{q+1}-}+\#
(C^{\lambda}_{T_{q+1}\log(1/{\lambda}
)}(X_{q+1}))
\simeq(1/{\lambda})^{1-z} +
(1/{\lambda})^z\simeq(1/{\lambda})^{\max(z,1-z)}$. We conclude that
${\tilde H}^{{\lambda},\mu}_{T_{q+1}}(X_{q+1})= \log(1+v^{{\lambda
},\mu}_{T_{q+1}})/\break \log(\#
((X_{q+1})_{{\lambda},\mu})) \simeq
\max(z,1-z)={\tilde H}_{T_{q+1}}(X_{q+1})$.
All of this can be done rigorously without difficulty
and we deduce that for ${\varepsilon}>0$ and all $\mu\in(0,1]$,
$\lim_{{\lambda}\to0}
{\mathbb{P}}_M(
|{\tilde H}^{{\lambda},\mu}_{T_{q+1}}(X_{q+1})-{\tilde
H}_{T_{q+1}}(X_{q+1})|\geq
{\varepsilon})=0$.
\end{itemize}

\textit{Case $z=1$.} Let $a,b \in{\mathcal B}_q \cup\{-A,A\}$ be such that
$D_{T_{q+1}-}(X_{q+1})=[a,b]$.
We assume that $a,b \in{\mathcal B}_q$, the other cases
being treated in a similar way.
We thus have $h_a:={\tilde H}_{T_{q+1}-}(a)>0$,
$h_b:={\tilde H}_{T_{q+1}-}(b)>0$.
We also have ${\tilde H}_{T_{q+1}}(x)={\tilde H}_{T_{q+1}-}(x)$
for all $x\in{\mathcal B}_q \setminus[a,b]$, ${\tilde H}_{T_{q+1}}(x)=1$
for all $x\in{\mathcal B}_q \cap(a,b)$, $Z_{T_{q+1}}(c)=Z_{T_{q+1}-}(c)$
for all $c\in{\mathcal C}_{q+1}$ with $c\cap(a,b)=\varnothing$ and
$Z_{T_{q+1}}(c)=0$ for all $c\in{\mathcal C}_{q+1}$ with $c\subset(a,b)$.

Let $\mu\in(0,1]$.
Now, consider ${\tilde\Omega}^{{\lambda},\mu}$, the event that
for all $c\in{\mathcal C}_{q}$ such that $c \subset(a,b)$, we have
${\tilde Z}^{{\lambda},\mu}_{T_{q+1}-}(c)=1$, that ${\tilde
H}^{{\lambda},\mu}_{T_{q+1}-}(a)>0$,
that ${\tilde H}^{{\lambda},\mu}_{T_{q+1}-}(b)>0$ and that
for all $x \in{\mathcal B}_q \cap(a,b)$, ${\tilde H}^{{\lambda},\mu
}_{T_{q+1}-}(x)=0$.
Then (a), (c) and Lemma \ref{undoncplein} collectively imply that
$\lim_{{\lambda}\to0} {\mathbb{P}}_M({\tilde\Omega}^{{\lambda
},\mu})=1$ for all $\mu\in(0,1]$.

\begin{itemize}
\item We can easily check that ${\mathcal E}^{{\lambda},\mu}_{q} \cap
{\tilde\Omega}^{{\lambda},\mu}
\subset
{\mathcal E}^{{\lambda},\mu}_{q+1}$
(because for $c \in{\mathcal C}_{q+1}$ with $c \subset[a$, $b]$,
we have $c_{{\lambda},\mu}\subset C^{\lambda}_{T_{q+1}\log
(1/{\lambda})-}(X_{q+1})$,
while for $c \in{\mathcal C}_{q+1}$ with $c \cap[a,b] = \varnothing$,
the vacant sites in $a_{{\lambda},\mu}$ and $b_{{\lambda},\mu}$
separate $c_{{\lambda},\mu}$ from\break
$ C^{\lambda}_{T_{q+1}\log(1/{\lambda})-}(X_{q+1})$). As a
consequence, $({\mathcal H}_{q+1})$(iii)
holds for all $\mu\in(0,1]$. 

\item On ${\tilde\Omega}^{{\lambda},\mu}$, we have ${\tilde
Z}^{{\lambda},\mu}
_{T_{q+1}}(c)=0=Z_{T_{q+1}}(c)$ for all
$c \in{\mathcal C}_{q+1}$ with $c \subset[a,b]$, and
${\tilde Z}^{{\lambda},\mu}_{T_{q+1}}(c)={\tilde Z}^{{\lambda},\mu
}_{T_{q+1}-}(c)$ for $c \in{\mathcal C}_{q+1}$
with $c\cap(a,b)=\varnothing$, from which $({\mathcal H}_{q+1})$(i)
easily follows
[using (a)].

\item We also have, still on ${\tilde\Omega}^{{\lambda},\mu}$, that
${\tilde H}^{{\lambda},\mu}_{T_{q+1}}(x)=1={\tilde H}_{T_{q+1}}(x)$
for all $x \in{\mathcal B}
_{q+1}$ with
$x \in(a,b)$, and $({\mathcal H}_{q+1})$(ii) follows for those $x$.
For $x \in{\mathcal B}_{q+1}$ with $x \notin[a,b]$, we have
${\tilde H}^{{\lambda},\mu}_{T_{q+1}}(x)={\tilde H}^{{\lambda},\mu
}_{T_{q+1}-}(x)$, hence $({\mathcal H}_{q+1})$(ii)
follows by point (c).

Finally, we have to check that $({\mathcal H}_{q+1})$(ii) holds for $x=a$
and $x=b$.
Consider, for example, the case of $a$.
Here, we are in the situation
where \mbox{$Z_{T_{q+1}}(a+)=0$} so that, of course, ${\tilde H}_{T_{q+1}}(a)
= 1$.
Let $c$ be the cell containing~$a+$.
We know that ${\tilde Z}^{{\lambda},\mu/2}_{T_{q+1}-}(c) =1$ which,
on ${\tilde\Omega}^{{\lambda},\mu/2}$,
implies that all sites between $a+\frac{\mu}{2\log(1/{\lambda})}$ and
$a+\frac
{\mu}{\log(1/{\lambda})}$, that is, on an interval of length $\frac
{\mu}{2
\log(1/{\lambda})}$, are empty at time $T_{q+1}$, showing that a fixed
proportion of $a_{{\lambda},\mu}$ is empty.
Recalling that $\lim_{{\lambda}\to0}{\mathbb{P}}_M( {\tilde\Omega
}^{{\lambda},\mu/2}) = 1$, it readily
follows that for all ${\varepsilon}>0$, $\lim_{{\lambda}\to
0}{\mathbb{P}}_M( {\tilde H}
^{{\lambda},\mu}
_{T_{q+1}}(a) > 1-{\varepsilon}) =1$.
Recalling that $ {\tilde H}^{{\lambda},\mu}_{T_{q+1}}(a) \leq1$, we
conclude that $({\mathcal H}
_{q+1})$(ii) holds for $x=a$.
\end{itemize}



\subsubsection*{Conclusion} Using points (b) and (e) above (with $q=0,\ldots,n$),
plus very similar arguments on the time interval $(T_n,T]$
(during which there are no fires), we deduce that
for all $x_0 \in(-A,A)\setminus{\mathcal B}_n$ and ${\varepsilon}>0$,
\[
\lim_{{\lambda}\to0}{\mathbb{P}}_M \biggl({\sup_{[0,T]}} |Z^{\lambda
}_t(x_0) - Z_t(x_0)| +
\int_0^T \delta(D^{\lambda}_t(x_0),D_t(x_0))\,dt \geq{\varepsilon}\biggr) =0.
\]
But, of course, for $x_0 \in(-A,A)$, we have ${\mathbb{P}}(x_0 \in
{\mathcal B}_n)=0$
so that
\[
\lim_{{\lambda}\to0}{\mathbb{P}}\biggl({\sup_{[0,T]}} |Z^{\lambda}_t(x_0)
- Z_t(x_0)| +
\int_0^T \delta(D^{\lambda}_t(x_0),D_t(x_0))\,dt \geq{\varepsilon}\biggr) =0.
\]
It remains to prove that for $t \in[0,T]$ and $x_0 \in(-A,A)$, we have
\[
\lim_{{\lambda}\to0} {\mathbb{P}}(\delta(D^{\lambda
}_t(x_0),D_t(x_0)))=0.
\]

\vspace*{3pt}

\textit{Case $t\ne1$.} We deduce from point (d) above that
if $x_0 \notin{\mathcal B}_n$ and $t \notin{\mathcal K}$, then we have
$\lim_{{\lambda}\to0} {\mathbb{P}}_M(\delta(D^{\lambda
}_t(x_0),D_t(x_0)))=0$. Since
${\mathbb{P}}(x_0 \in{\mathcal B}_n)=0$ and ${\mathbb{P}}(t \in
{\mathcal K})=0$ (because $t \ne1$,
recalling the definition of ${\mathcal K}$), we easily arrive at the
desired conclusion.\vspace*{8pt}

\textit{Case $t=1$.} In this case, $t \in{\mathcal K}$, but the result
still holds.
Observe that \mbox{$Z_1(x_0)=1$}, by construction.
Consider $q\in\{0,\ldots,n\}$ such that $T_q<1<T_{q+1}$ (with the
convention that
$T_0=0$, $T_{n+1}=T$) and consider $a,b \in{\mathcal B}_q \cup\{-A,A\}$
such that $D_1(x_0)=[a,b]$.
Using the same arguments as in the proof of (d) (see Step 1),
we then easily check
that $\lim_{{\lambda}\to0} {\mathbb{P}}_M( D^{\lambda}_{1}(x_0)
\subset[a-{\varepsilon
},b+{\varepsilon}])=1$
for all ${\varepsilon}>0$ (the set ${\mathcal K}$ was not considered there).
We also check, as in the proof of (d) (see Step 2), that for all $y \in
{\mathcal B}_q$
with $y \in(a,b)$,
$\lim_{{\lambda}\to0} {\mathbb{P}}_M( H^{{\lambda},1}_{1}(y)=0)=1$
[the set under
consideration there was ${\mathcal K}$, but the time $1$ was not useful since
$1$ is a.s. not a time where some $H(x)$ reaches $0$ for the first time].
Finally, we just have to prove that for all $c \in{\mathcal C}_q$ with
$c \subset(a,b)$, $\lim_{{\lambda}\to0} {\mathbb{P}}_M({\tilde
Z}^{{\lambda},1}_1(c)=1)=1$.
Thus, let $c \in{\mathcal C}_q$ with $c \subset(a,b)$ and
recall that $\lim_{{\lambda}\to0}{\mathbb{P}}_M({\mathcal
E}^{{\lambda},1}_q)=1$. However, on
${\mathcal E}^{{\lambda},1}_q$, there are no death events in
$c_{\lambda}$ during the time interval
$[0,\log(1/{\lambda})]$, so each site of $c_{{\lambda},1}$ is
occupied at time
$\log(1/{\lambda})$ with probability $1- {\lambda}$ and, hence, all
the sites
of $c_{{\lambda},1}$ are occupied with probability $(1-{\lambda})^{\#
(c_{{\lambda},1})}$.
Since $\#(c_{{\lambda},1}) \leq2A/({\lambda}\log(1/{\lambda}))$,
we get
${\mathbb{P}}_M({\tilde Z}^{{\lambda},1}_1(c)=1 \vert{\mathcal
E}^{{\lambda},1}_q) \geq(1-{\lambda})^{2A/({\lambda}
\log
(1/{\lambda}))}$,
which tends to $1$ as ${\lambda}$ tends to $0$.
Since we know that $\lim_{{\lambda}\to0}{\mathbb{P}}_M({\mathcal
E}^{{\lambda},1}_q)=1$,
we deduce that $\lim_{{\lambda}\to0} {\mathbb{P}}_M( [a+1/\log
(1/{\lambda}),
b-1/\log(1/{\lambda})] \subset D^{\lambda}_{1}(x_0))=1$.

Finally,
$\lim_{{\lambda}\to0} {\mathbb{P}}_M( \delta(D^{\lambda
}_{1}(x_0),D_1(x_0))\geq
{\varepsilon})=0$
for all ${\varepsilon}>0$, which was our goal.
\end{pf*}

\section{Cluster size distribution}\label{conseq}

The aim of this section is to prove Corollary \ref{coco}. We
will use Theorem \ref{converge}, which asserts that the ${\lambda}$-FFP
behaves like the LFFP for ${\lambda}>0$ small enough. We start
with preliminary results.
\begin{lem}\label{zunif}
Consider an LFFP $(Z_t(x),D_t(x),H_t(x))_{t\geq0, x\in{\mathbb{R}}}$.
We then have the following:

{\smallskipamount=0pt
\begin{longlist}
\item
for any $t\in(1,\infty)$, $x\in{\mathbb{R}}$ and
$z\in[0,1)$, ${\mathbb{P}}[Z_t(x)=z]=0$;

\item for any $t\in[0,\infty)$, $B>0$ and $x\in{\mathbb{R}}$,
$P[|D_t(x)|=B]=0$;

\item there are constants $C>0$ and $\kappa_1>0$ such that for all
$t\in[0,\infty)$, $x\in{\mathbb{R}}$ and $B>0$,
${\mathbb{P}}[|D_t(x)|\geq B] \leq C e^{-\kappa_1 B}$;

\item there are constants $c>0$ and $\kappa_2>0$ such that for all
$t\in[3/2,\infty)$, $x\in{\mathbb{R}}$ and $B>0$,
${\mathbb{P}}[|D_t(x)|\geq B] \geq c e^{-\kappa_2 B}$;

\item there exist constants $0<c<C$ such that for all $t\geq5/2$,
$0\leq a < b < 1$ and $x\in{\mathbb{R}}$,
$c(b-a) \leq{\mathbb{P}}(Z_t(x)\in[a,b]) \leq C(b-a)$.
\end{longlist}}
\end{lem}
\begin{pf} By translation invariance,
it suffices to treat the case $x=0$.\vspace*{8pt}

\textit{Point} (i). By Definition \ref{dflffp},
we see that for $t\in[0,1]$,
we have a.s. $Z_t(0)=t$. However, for $t> 1$ and $z\in[0,1)$,
$Z_t(0)=z$ implies that the cluster containing $0$
has been killed at time $t-z$, so, necessarily,
$M(\{t-z\}\times{\mathbb{R}})>0$. This happens with probability $0$
since $t-z$ is deterministic.\vspace*{8pt}

\textit{Point} (ii).
Recalling Definition \ref{dflffp}, we see that for any $t\in[0,T]$,
$|D_t(0)|$ is either $0$ or of the form $|X_i-X_j|$ (with $i\ne j$),
where $(T_i,X_i)_{i\geq1}$ are the marks of the
Poisson measure $M$. As before, we easily conclude that for $B>0$,
${\mathbb{P}}(|D_t(0)|=B)=0$.\vspace*{8pt}

\textit{Point} (iii).
First, if $t \in[0,1)$, then we have a.s. $|D_t(0)|=0$ and the result
is obvious.
Next, consider
$t\geq1$. Recalling Definition \ref{dflffp}, we see that $|D_t(0)|= |L_t(0)|
+R_t(0)$. Clearly, $|L_t(0)|$ and $R_t(0)$ have the same law.
For $B>0$, $\{R_t(0)> B \} \subset
\{ M([t-1/4,t]\times[0,B])=0\}$.
Indeed, on $\{ M([t-1/4,t]\times[0,B])>0\}$,
denote by $(\tau,X)\in[t-1/4,t]\times[0,B]$ a mark of $M$. Then, either:

$\bullet$
$Z_{\tau-}(X)=1$, in which case this mark starts a macroscopic fire so that
$Z_{\tau}(X)=0$ and $Z_{s}(X)=s-\tau<1$ for all $s\in[\tau,\tau+1)$
(since $\tau\in[t-1/4,t]$, we clearly have $t \in[\tau,\tau+1)$
so that $Z_t(X)<1$ and, as a consequence, $R_t(0)\leq X \leq B$); or

$\bullet$
$Z_{\tau-}(X) \in(1/4,1]$ so that $H_{\tau}(X)=Z_{\tau-}(X)$
and thus $H_{s}(X)=Z_{\tau-}(X) - (s-\tau) >0$ for all $s\in[\tau,
\tau+
Z_{\tau-}(X))$ (since $\tau\in[t-1/4,t]$ and $Z_{\tau-}(X)>1/4$,
we have $t \in[\tau, \tau+Z_{\tau-}(X))$,
so $H_t(X)>0$ and, hence, $R_t(0)\leq X \leq B$); or, finally,

$\bullet$
$Z_{\tau-}(X)\leq1/4$, in which case
$Z_{s}(X)=Z_{\tau-}(X)+(s-\tau)<1$
for all $s\in(\tau, \tau+1-Z_{\tau-}(X))$ and,
in particular, $Z_t(X)<1$, hence $R_t(0)\leq X \leq B$.

As a conclusion, for all $t\geq1$, ${\mathbb{P}}[R_t(0) > B] \leq
{\mathbb{P}}[ M([t-1/4,t]\times[0,B])=0]=e^{-B/4}$,
so ${\mathbb{P}}[|D_t(0)|> B] \leq
{\mathbb{P}}[|L_t(0)|> B/2]+ {\mathbb{P}}[R_t(0) > B/2] \leq2 e^{-B/8}$.\vspace*{8pt}

\textit{Point} (iv). We first observe that for all $(t_0,x_0)$ such that
$M(\{t_0,x_0\})=1$,
we have $\max(1-Z_{t}(x_0),H_{t}(x_0))>0$ for all $t\in[t_0,t_0+1/2)$.

Indeed, if $Z_{t_0-}(x_0)=1$, then $Z_{t_0+s}(x_0) \leq
s<1$ for all $s\in[0,1)$. If, now, $z=Z_{t_0-}(x_0)<1$, then
$Z_{t_0+s}(x_0) = s+z<1$ for $s \in[0,1-z)$ and $H_{t_0+s}(x_0) = z-s>0$
for $s\in[0,z)$ so that
$\max(1-Z_{t_0+s}(x_0),H_{t_0+s}(x_0))>0$ for all $s\in[0,1/2)$.

Once this is seen, fix $t\geq3/2$. Consider the event
${\tilde\Omega}_{t,B}={\tilde\Omega}^1_{t,B} \cap{\tilde\Omega
}^2_{t}\cap{\tilde\Omega}^3_{t,B}$,
where:

$\bullet$ ${\tilde\Omega}^{1}_{t,B}=\{M([t-3/2,t]\times[0,B])=0\}$;

$\bullet$ ${\tilde\Omega}^2_{t}$ is the event that in the box
$[t-3/2,t]\times[-1,0]$,
$M$ has exactly four marks, $(S_i,Y_i)_{i=1,\ldots,4}$, with
$Y_4<Y_3<Y_2<Y_1$, $t-3/2 < S_1 < t-1$, $S_1<S_2 < S_1+1/2$, $S_2<S_3< S_2+1/2$,
$S_3<S_4< S_3+1/2$ and $S_4+1/2>t$.

$\bullet$ ${\tilde\Omega}^3_{t,B}$ is the event that in the box
$[t-3/2,t]\times[B,B+1]$,
$M$ has exactly four marks, $({\tilde S}_i,{\tilde Y}_i)_{i=1,\ldots,4}$, with
${\tilde Y}_1<{\tilde Y}_2<{\tilde Y}_3<{\tilde Y}_4$, $t-3/2 < {\tilde
S}_1 < t-1$, ${\tilde S}_1<{\tilde S}_2 <
{\tilde S}
_1+1/2$, ${\tilde S}_2<{\tilde S}_3< S_2+1/2$,
${\tilde S}_3<{\tilde S}_4< {\tilde S}_3+1/2$ and ${\tilde S}_4+1/2>t$.

Of course, we have $p:={\mathbb{P}}({\tilde\Omega}^2_{t})={\mathbb
{P}}({\tilde\Omega}^3_{t,B})>0$
and this probability does not depend on $t\geq3/2$ or on $B>0$.
Furthermore, ${\mathbb{P}}({\tilde\Omega}^1_{t,B})=e^{-3B/2}$. These three
events being independent, we conclude that
${\mathbb{P}}({\tilde\Omega}_{t,B}) \geq p^2 e^{-3B/2}$. To conclude
the proof of (iv),
it thus suffices to check that ${\tilde\Omega}_{t,B}\subset\{
[0,B]\subset
D_t(0)\}$.
However, on ${\tilde\Omega}_{t,B}$, using the arguments described at
the beginning
of the proof of point (iv), we observe that:

$\bullet$ the fire starting at $(S_2,Y_2)$ cannot affect $[0,B]$
because at time $S_2 \in[S_1, S_1 +1/2)$, $H_{S_2}(Y_1)>0$ or
$Z_{S_2}(Y_1)>0$, with $Y_2<Y_1<0$;

$\bullet$ then the fire starting at $(S_3,Y_3)$ cannot affect $[0,B]$ because
at time $S_3 \in[S_2,S_2+1/2)$, $H_{S_3}(Y_2)>0$ or
$Z_{S_3}(Y_2)>0$, with $Y_3<Y_2<0$;

$\bullet$ then the fire starting at $(S_4,Y_4)$ cannot affect $[0,B]$ because
at time $S_4 \in[S_3,S_3+1/2)$, $H_{S_4}(Y_3)>0$ or
$Z_{S_4}(Y_3)>0$, with $Y_4<Y_3<0$;

$\bullet$ furthermore, the fires starting to the left of $-1$ during $(S_1,t]$
cannot affect $[0,B]$ because for all $t\in(S_1,t]$, there is
always a site $x_t \in\{Y_1,Y_2,Y_3,Y_4\} \subset[-1,0]$
with $H_t(x_t)>0$ or $Z_t(x_t)<1$;

$\bullet$ the same arguments apply on the right of $B$.

As a conclusion, the zone $[0,B]$ is not affected by any fire
during $(S_1 \lor{\tilde S}_1,t]$. Since the length of this time
interval is
greater than $1$, we deduce that for all $x \in[0,B]$,
$Z_t(x)=\min(Z_{S_1 \lor{\tilde S}_1} + t- S_1 \lor{\tilde S}_1,1)
\geq\min(t- S_1 \lor{\tilde S}_1,1)=1$
and $H_t(x)=\max(H_{S_1 \lor{\tilde S}_1} - (t- S_1 \lor{\tilde S}_1),0)
\leq\max(1 - (t- S_1 \lor{\tilde S}_1),0)=0$,\vspace*{2pt} hence that
$[0,B]\subset D_t(0)$.\vspace*{8pt}

\textit{Point} (v). We observe, recalling Definition \ref{dflffp},
that for $0\leq a < b < 1$ and $t\geq1$,
we have $Z_t(0)\in[a,b]$ if and only there exists $\tau\in[t-b,t-a]$
such that $Z_\tau(0)=0$. This happens if and only if
$X_{t,a,b}:=\int_{t-b}^{t-a}\int_{\mathbb{R}}{\mathbf{1}}_{\{y\in
D_{{s-}}(0)\}
}M(ds,dy)\geq1$.
We deduce that
\[
{\mathbb{P}}\bigl(Z_t(0)\in[a,b] \bigr)
={\mathbb{P}}(X_{t,a,b}\geq1 )
\leq{\mathbb{E}}[X_{t,a,b} ]
= \int_{t-b}^{t-a}{\mathbb{E}}[|D_s(0)|]\,ds \leq C(b-a),
\]
where we have used point (iii) for the last inequality.

Next, we have $\{M([t-b,t-a]\times D_{t-b}(0))\geq1 \}\subset
\{X_{t,a,b}\geq1\}$: it suffices to note that a.s.
$\{X_{t,a,b}=0\} \subset\{X_{t,a,b}=0, D_{t-b}(0)\subset D_s(0)$ for
all $s\in
[t-b,t-a] \} \subset\{M([t-b,t-a]\times D_{t-b}(0))=0\}$.
Now, since $D_{t-b}(0)$ is ${\mathcal F}^M_{t-b}$-measurable, we deduce
that for $t\geq5/2$,
\begin{eqnarray*}
{\mathbb{P}}\bigl(Z_t(0)\in[a,b] \bigr)
&\geq&{\mathbb{P}}\bigl[M\bigl((t-b,t-a]\times D_{t-b}(0)\bigr)>0 \bigr]\\
&\geq&{\mathbb{P}}[ |D_{t-b}(0)|\geq1 ] \bigl(1-e^{-(b-a)}\bigr) \geq c \bigl(1-e^{-(b-a)}\bigr),
\end{eqnarray*}
where we have used point (iv) (here, $t-b\geq3/2$) to get the last inequality.
This completes the proof since $1-e^{-x}\geq x/2$ for all $x\in[0,1]$.
\end{pf}

We now may tackle the following proof.
\begin{pf*}{Proof of Corollary \ref{coco}}
We thus consider, for each ${\lambda}>0$, a ${\lambda}$-FFP $(\eta
^{\lambda}
_t)_{t\geq0}$.
Also, let $(Z_t(x),D_t(x),H_t(x))_{t\geq0, x\in{\mathbb{R}}}$ be an LFFP.\vspace*{8pt}

\textit{Point} (i). Using Lemma \ref{zunif}(v),
we only need to prove that
for all $0\leq a < b < 1$ and all $t\geq5/2$,
\[
\lim_{{\lambda}\to0} {\mathbb{P}}\bigl(\#\bigl(C^{\lambda}_{t\log
(1/{\lambda})}(0)\bigr)\in[{\lambda}^{-a},{\lambda}^{-b}]
\bigr)
={\mathbb{P}}\bigl(Z_t(0)\in[a,b] \bigr).
\]
Recalling (\ref{zlambda}), we observe that
\[
{\mathbb{P}}\bigl(\#\bigl(C^{\lambda}_{t\log(1/{\lambda})}(0)\bigr) \in[{\lambda
}^{-a},{\lambda}^{-b}] \bigr) =
{\mathbb{P}}\bigl(Z^{\lambda}_t(0) \in[a+{\varepsilon}(a,{\lambda
}),b+{\varepsilon}(b,{\lambda}
)] \bigr),
\]
where ${\varepsilon}(z,{\lambda})= \log(1+{\lambda}^z)/\log
(1/{\lambda})
\to0$ as ${\lambda}\to0$ (if $z\geq0$).

We arrive at the desired conclusion by using Theorem \ref{converge}
[which asserts that $Z_t^{\lambda}(0)$ goes in law to
$Z_t(0)$] and Lemma \ref{zunif}(i) [from which
${\mathbb{P}}(Z_t(0)=a)={\mathbb{P}}(Z_t(0)=b)=0 $].\vspace*{8pt}

\textit{Point} (ii). Using part (iv) of Lemma \ref{zunif}(iii)
and recalling (\ref{dlambda}), it suffices
to check that for all $t\geq3/2$ and all $B>0$, we have
\[
\lim_{{\lambda}\to0} {\mathbb{P}}[|D^{\lambda}_t(0)|\geq B ] =
{\mathbb{P}}[|D_t(0)|\geq B ].
\]
This follows from Theorem \ref{converge} and the fact that
${\mathbb{P}}(|D_t(0)|=B)=0$, thanks to Lemma \ref{zunif}(ii).
\end{pf*}

\section*{Acknowledgment}
We are grateful to the referee who helped us to make the proofs more
readable and, indeed, correct.


%
\printaddresses

\end{document}